\documentclass[11pt,reqno]{amsart}

 \usepackage{curves} 
\usepackage{graphicx} 
\usepackage{amsmath, amsfonts, amssymb, amscd,graphics, graphpap, tikz, color, xcolor, cancel,enumerate,enumitem, hyperref, mathabx, mathrsfs, mathtools}
\usepackage[scr=rsfs,cal=boondox]{mathalfa}

\usetikzlibrary{arrows.meta,arrows}
\usetikzlibrary{decorations.pathreplacing,decorations.markings}
\tikzset{arrow data/.style 2 args={%
      decoration={%
         markings,
         mark=at position #1 with \arrow{#2}},
         postaction=decorate}
      }%
      
\usepackage[mathscr]{euscript}
  \usepackage[normalem]{ulem}
  \usepackage{soul}

\usepackage{curves} 
\usepackage{graphicx}

\makeatletter
\newcommand{\doublewidetilde}[1]{{%
  \mathpalette\double@widetilde{#1}%
}}

\textwidth= 
148mm
\numberwithin{equation}{section}

\theoremstyle{plain}
\newtheorem{theo}{Theorem}[section]
\newtheorem{lem}[theo]{Lemma}
\newtheorem{prop}[theo]{Proposition}

\newtheorem{cor}[theo]{Corollary}

\theoremstyle{definition}

\newtheorem{example}[theo]{Example}
\newtheorem{definition}[theo]{Definition}

\newenvironment{pf}{\noindent{\it Proof.\,}}{\hfill $\square$\par \medskip}

\theoremstyle{plain}

\theoremstyle{definition}

\newcommand{\rank}{\operatorname{rank}}
\newcommand{\beq}{\begin{equation}}
\newcommand{\eeq}{\end{equation}}
\renewcommand{\a}{\alpha}
\renewcommand{\b}{\beta}

\renewcommand{\d}{\delta}

\newcommand{\f}{\varphi}
\newcommand{\g}{\gamma}
\newcommand{\h}{\eta}

\renewcommand{\o}{\omega}

\newcommand{\s}{\sigma}
\renewcommand{\t}{\tau}

\newcommand{\bR}{\mathbb{R}}

\newcommand{\bA}{\mathbb{A}}

\newcommand{\bT}{\mathbb{T}}
\newcommand{\bV}{\mathbb{V}}

\newcommand\U{\mathrm{U}}

\newcommand{\cC}{\mathcal{C}}
\newcommand{\cD}{\mathscr{D}}
\newcommand{\cE}{\mathscr{E}}

\newcommand{\cG}{\mathscr{G}}

\newcommand{\cK}{\mathscr{K}}
\newcommand{\cL}{\mathscr{L}}
\newcommand{\cM}{\mathscr{M}}

\newcommand{\cQ}{\mathscr{Q}}
\newcommand{\cR}{\mathscr{R}}
\newcommand{\cS}{\mathscr{S}}

\newcommand{\cU}{\mathscr{U}}
\newcommand{\cV}{\mathscr{V}}

\newcommand{\p}{\partial}

\renewcommand{\square}{\kern1pt\vbox
{\hrule height 0.6pt\hbox{\vrule width 0.6pt\hskip 3pt
\vbox{\vskip 6pt}\hskip 3pt\vrule width 0.6pt}\hrule height0.6pt}\kern1pt}

\newcommand{\wt}{\widetilde}

\newcommand{\wh}{\widehat}
\newcommand{\wc}{\widecheck}

\newcommand{\bt}{\begin{theo}\ \ }
\newcommand{\et}{\end{theo}}
\newcommand{\bp}{\begin{prop}\ \ }
\newcommand{\ep}{\end{prop}}
\newcommand{\bc}{\begin{cor}\ \ }
\newcommand{\ec}{\end{cor}}
\newcommand{\bl}{\begin{lem}\ \ }
\newcommand{\el}{\end{lem}}
\newcommand{\bd}{\begin{definition}}
\newcommand{\ed}{\end{definition}}
\newcommand{\be}{\begin{equation}}
\newcommand{\ee}{\end{equation}}

\def\<#1,#2>{\langle\,#1,\,#2\,\rangle}
\newcommand{\arr}{\begin{array}{rlll}}
\newcommand{\ea}{\end{array}}
\newcommand{\bea}{\begin{eqnarray}}
\newcommand{\eea}{\end{eqnarray}}
\newcommand{\bean}{\begin{eqnarray*}}
\newcommand{\eean}{\end{eqnarray*}}

\renewcommand{\=}{:=}

\newcommand{\ve}{\varepsilon}

  \newcommand{\vertiii}[1]{{\left\vert\kern-0.25ex\left\vert\kern-0.25ex\left\vert #1 
    \right\vert\kern-0.25ex\right\vert\kern-0.25ex\right\vert}}

\newcommand{\Attain}{\cM\text{\it -Att}}

\newcommand{\Reach}{\cR\text{\it each}}

\def\sideremark#1{\ifvmode\leavevmode\fi\vadjust{
\vbox to0pt{\hbox to 0pt{\hskip\hsize\hskip1em
\vbox{\hsize3cm\tiny\raggedright\pretolerance10000
\noindent #1\hfill}\hss}\vbox to8pt{\vfil}\vss}}}

\title[Distributions and controllability  problems (II)]{Distributions and controllability  problems (II)}
 \author[Cristina Giannotti,  Andrea Spiro and Marta Zoppello]{Cristina Giannotti \quad  Andrea Spiro \quad Marta Zoppello}
  \thanks{{\it Data availability statement.} 
No datasets were generated or analysed during the current study.}

\begin{document}

\begin{abstract} 
In [C. Giannotti, A. Spiro, M. Zoppello, {\it Distributions and controllability  problems (I)}, preprint posted on ArXiv  2401.07560 (2024)], we introduced a new  approach to the real analytic  non-linear  control systems of the form $\dot q^i = f^i(t, q, w)$, with  controls $w = (w^\a)$ running in a connected open set  $\cK$ of  $ \bR^m$  and   states represented by points  $q = (q^i)$ in a configuration space  $\cQ \= \bR^n$.  The new  approach   consists of a  differential-geometric study of  (a)  the oriented piecewise regular curves in the {\it extended   space-time}  $\cM = \bR \times \cQ \times \cK$,  which are the (completed) graphs  of the  piecewise real analytic solutions $t \mapsto (q(t), w(t))$ of the control system,  and  (b)  the local structure of the sets  of points of $\cM$ that are reachable from an initial point  $x_o = (t_o, q_o, w_o)  \in \cM$ through such   (completed) graphs. The main results of that paper   are two new criterions which can be used to  establish the  small time local controllability near stable points of real analytic non-linear systems.  The goal of this paper is to offer a friendly user's guide to those   criterions,   illustrating  them by  several examples. In particular, we  analyse  certain  non-linear control systems, for which the new criterions    show   that  they are   small time locally  controllable at their  stable points, while, at the best  of our knowledge,   all other previous criterions are either inconclusive or  not applicable. 
\end{abstract}
\maketitle
\setcounter{section}{0}
\section{Introduction}
This paper is the natural continuation of \cite{GSZ2}, where we  investigated  the  accessibility  and the  small time controllability of    real analytic systems  $\dot q = f(t, q, w)$ with controls $w = (w^\a)$ in an open subset $\cK$ of $\bR^m$ and controlled states $q = (q^i)$ in an affine space $\cQ = \bR^n$.  There we considered completed graphs   in $\cM = \bR \times \cQ \times \cK \subset \bR^{1 + n + m}$  of  solutions $t \mapsto (q(t), w(t))$ with piecewise real analytic controls  (i.e. graphs of  solutions to which we add  arcs with trivial projections on $\bR \times \cQ$ just    to make them connected) and we characterised them in purely differential geometric terms. Indeed, such completed graphs coincide with the piecewise real analytic oriented  curves $\h(t)$ which can be obtained by  composing a finite number of  smooth arcs with  tangent vectors $\dot \h(t)$
in  the distribution $\cD$  generated  by the vector fields 
$$\bT^o = \frac{\p}{\p t} + f^i(t, q, w) \frac{\p}{\p q^i}\ ,\qquad W^o_1 = \frac{\p}{\p w^1} \ ,\ \ \ldots \ \  , \ \ W^o_m = \frac{\p}{\p w^m}\ ,$$
  with either nowhere vanishing projections onto the quotients $\cD_{\h(t)}\!\!\! \mod \langle W_\a^o\rangle_{1 \leq \a \leq m}$ or with  identically vanishing projections on such quotients. The points in $\cM$, which can be reached from  a given $x_o = (t_o, q_o, w_o) $  through  oriented curves of this kind, are  called  {\it $\cM$-attainable points from $x_o$}. The first main result of  \cite{GSZ2} consists in the proof that, given an $\cM$-attainable  point $y_o$   from  a fixed $x_o$, there exists  a  large set of  points, which are  not only  $\cM$-attainable  from the same $x_o$, but also reachable  directly from $y_o$    by applying  an appropriate  composition of flows of certain vector fields, named    {\it surrogate   fields}. These  vector   fields   take values in   the so-called  {\it secondary distribution} $\cD^{II} \subset T \cM$, that is the (possibly non-regular) distribution which is generated by the real analytic vector fields 
$$\underset{\text{$k$-times}}{\underbrace{[\bT^o, [\bT^o, \ldots, [ \bT^o, }}W^o_\b]\ldots]]\ , \qquad k = 0, 1, 2, 3, \ldots\ .$$
This result  implies that for each $\cM$-attainable point $y_o$ from $x_o$, there exists  a large set of other $\cM$-attainable points from the same $x_o$, which is included in the  orbit of $y_o$ under the action of the pseudogroup of local diffeomorphisms $\cG^{II}$  generated by the flows of the vector fields in $\cD^{II}$.  The second main result of \cite{GSZ2} consists in a couple of sufficient conditions, under which all points of a  sufficiently small local $\cG^{II}$-orbit of an $\cM$-attainable point $y_o$  from $x_o$ are also  reachable through   flows of surrogate vector fields  and thus  $\cM$-attainable from  $x_o$ as well. Since the  local $\cG^{II}$-orbits are in many cases  easily computable  via the  Chow-Rashevski\u\i\ and   Frobenius Theorems, whenever one of such sufficient conditions is satisfied, the local behaviour   of the $\cM$-attainable points and of their projections on   $\cQ = \bR^n$ (= the  reachable states through piecewise real analytic controls) can  be  determined with little effort. \par
\smallskip
The purpose of this paper is to illustrate how this circle of ideas  works through explicit examples  and to indicate in which way they can be used to  establish  accessibility  and   small time local controllability in proximity of  stable states. \par
\smallskip
With this aim in mind,  we first  review the main notions that are necessary to state in detail  the   results of \cite{GSZ2}.  On this regard, we  stress the fact that in \cite{GSZ2} 
the discussion was designed to give  a purely geometric presentation   of all involved ingredients,  with the purpose of paving the way to  developments  in  much more  general settings. In particular,    all notions and properties are   there presented using   a   coordinate-free language. Here,  since we have a different goal, i.e.  
 providing a  friendly user's guide for   a widest as possible audience,  and  since  we are  working   just on open subsets  of  appropriate affine spaces  $\bR^N$,  whenever possible  we  re-formulate the contents of \cite{GSZ2}  using  just representations in    cartesian coordinates.
 Second,  we  very briefly recall   the Chow-Rashevski\u\i\ Theorem  (following  Rashevski\u\i's original  paper \cite{Ra, GSZ1})  and the  two main results in \cite{GSZ2}  that allow to  relate the orbits of the pseudo-group $\cG^{II}$ with the sets of  $\cM$-attainable points. \par
\smallskip
Third, we get into the actual core of this paper, namely a discussion   of several  applications of our main results. In particular, we  offer a new direct proof of  the classical Kalman Theorem on linear control systems  (see e.g. \cite{Bl, Co, Ju} and references therein)   and      prove the    small-time locally  controllability near  stable states for a particular non-linear control system, for which, at the best  of our knowledge,   all previously known criterions are either inconclusive or  not applicable. We conclude  proving that the   nonholonomic control systems given by  the    {\it controlled Chaplygin sleigh} (\cite{BKMM, BK}) and its {\it hydrodynamical variant}  (\cite{Ba, FG,  FGV, SZ}),  are both  small-time locally  controllable near their  states of equilibrium, a result which, despite of the simplicity of the  two systems, seems to us that it was till now   missing. \par
\smallskip
The structure of the paper is the following.  The two initial sections \S \ref{sect2} and \ref{sect3} are devoted to recall the  definitions  of  $\cM$-attainable points, their relation with the reachable points of  a control system and   an explanation of what  are the  secondary distribution and  the surrogate fields of a real analytic control system. As we mentioned above, most of the results of  \cite{GSZ2} are here re-presented using   standard  cartesian coordinates of an affine space.  Section \S \ref{sect3} ends with the definition of a very important class of  sets  of  $\cM$-attainable points for a given $x_o \in \cM$,  the so-called {\it  surrogate leaflets}.  In section \S \ref{sect4}, we recall  Rashevski\u\i's Theorem and  review  the basic relations between a set of  $\cM$-attainable points, the surrogate leaflets centred at its points and the local orbits of the same  points under the action of  the pseudogroup $\cG^{II}$, which is generated by the local diffeomorphisms of the flows of the vector fields in the secondary distribution.  At the end of that  section, we  review  the two  main   results in \cite{GSZ2},  which  allow to reduce the study of the surrogate leaflets to the analysis of the local orbits of $\cG^{II}$. 
The advertised applications are  in \S \ref{applications}.
\par
\medskip

\section{Control  systems, reachable sets and $\cM$-attainable sets} \label{thirdsection}
\label{sect2}
\subsection{Graphs of solutions of a control system and minimally completed graphs}
In this paper we consider    first order systems of control  equations on curves  $q(t) = (q^i(t))$  in  a   space $\cQ$ of  {\it states} of  a controlled dynamical system. We assume  that (a)  $\cQ$  is  an affine space  $\cQ = \bR^n$,  (b)  the   system  of differential equations that controls the dynamical  system is  real analytic and in normal form
\beq \label{controlsystem} 
\dot q^i(t) =  f^i(t, q(t), w(t))\ , \qquad 1 \leq i \leq n\ ,\eeq
and (c)     the  controls  $w(t) = (w^\a(t))_{ 1 \leq \a \leq m}$ are represented by   curves in an  open  connected region $\cK$  of $\bR^m$.    The {\it solutions}   of  \eqref{controlsystem}  are the maps  $t \mapsto (q(t), w(t)) \in \cQ \times \cK$, $t \in [a, b] \subset \bR$, in which  the  {\it control curve} $t \mapsto w(t)$  is a   measurable  map,  and the {\it states curve}   $t \mapsto q(t)$   is an absolutely continuous   map with the property that  the map  $t \mapsto (t, q(t), w(t))$  satisfies \eqref{controlsystem} at almost all points.  As in \cite{GSZ2}, we limit out discussion only to the  solutions, in which $w(t)$ is   piecewise $\cC^\o$, 
 possibly not continuous but     with  only finite  jumps  at the points of discontinuity, a condition that we   tacitly assume throughout the paper.\par
 \smallskip
 The {\it   extended space-time} of the   system \eqref{controlsystem}  is the cartesian product $\cM \= \bR \times \cQ \times \cK$ and the  {\it  graph  of a  solution} $t \mapsto \big( q(t), w(t)\big)$, $t \in [a,b]$,     is the subset  of $\cM$ given by the triples 
 $ \big(t, q(t), w(t)\big)$, $t \in [a,b]$. 
 By  the above definition of ``solution'',   any such  graph  is    a   (possibly not connected)  finite union of oriented  connected curves $\g_1, \ldots, \g_r$,  where each map $\g_j: I_j \subset \bR \to \cM$ is defined on an interval    $I_j  \subset \bR$  with possibly either one or both endpoints removed  and satisfies   the following  two conditions: (1) for each $\g_j$, $1 \leq j \leq r$,  the projection onto $\bR$  coincides with the interval $I_j$;  (2) the orientation  coincide with the   one  given by the  natural parameterisation $t \mapsto (t, q(t),  w(t))$. \par
 \smallskip
   A {\it  completed graph} of a solution is a  connected piecewise regular  oriented curve  $\h_1 \ast \h_2 \ast  \ldots\ast \h_{2r+1}$ in  $\bR \times \cQ \times  \cK$,  in which 
 \begin{itemize}[leftmargin = 20 pt]
 \item[(a)]  The union  of the {\it odd}  arcs $\h_{2 i +1}$, $ 0 \leq i \leq r$, 
  is the graph of a solution; 
  \item[(b)] Each {\it even } arc $\h_{2 i}$, $ 1 \leq i \leq r$,  is  a piecewise real analytic curve  with trivial projection onto $\bR \times \cQ$, which  joins the    endpoints (or  limit endpoints)  $(t_o, q_o, w_o)$,   $(t_o, q_o, w_o')$ 
of the two adjacent arcs $\h_{2 i - 1}$ and $\h_{2 i +1}$. 
 \end{itemize}
In other words, a completed graph is the union of the connected components of the graph of a solution to which are added some new arcs,  trivially projecting  on the time-axis $\bR$ and the states space $\cQ$, with the only purpose of making the graph a connected set   (see Fig. 1 and Fig. 2). 
 \par
 \smallskip
 \centerline{
\begin{tikzpicture}
  \begin{scope}[white]
        \draw[fill=red!30, semitransparent] (0.5, 0) -- (0.5,2.5) -- (1,4) --  (1,1.5) --cycle;
        \draw[fill=red!30, semitransparent] (2, 0) -- (2,2.5) -- (2.5,4) --  (2.5,1.5) --cycle;
        \draw[fill=red!30, semitransparent] (3.5, 0) -- (3.5,2.5) -- (4,4) --  (4,1.5) --cycle;
\end{scope}
\draw[->, line width = 0.6, blue]  (0,0.5) to (5,0.5);
\draw[->, line width = 0.6, red]  (0.5,0) to (0.5,4);
\draw[->, line width = 0.2, red]  (0.4,0.2) to (1.2,2.6);
\node at  (5, 0.3) { \color{blue} \tiny$q$};
\node at  (0.3, 3.75) { \color{red} \tiny $w^1$};
\node at  (1.5, 2.5) { \color{red}  \tiny $ w^2$};
\draw [line width = 0.8, red, dashed](0.6,0.8)  to  [out=40, in=250] (0.9, 1.9)   ; 
\draw [line width = 0.8, red, dashed](0.9,2.2)  to  [out=87, in=270] (0.95, 3.8)   ; 
\draw[fill, red]  (0.9, 1.9) circle [radius = 0.04];
\draw[fill, red]  (0.9, 2.2) circle [radius = 0.04];
\draw [line width = 0.7, blue](0.6,0.8)  to  [out=-5, in=230] (2.2, 1.9)   ; 
\draw [line width = 0.7, blue](2.22,2.2)  to  [out=60, in=240] (4, 3.8)   ; 
\draw[fill, blue]  (2.2, 1.9) circle [radius = 0.04];
\draw[fill, blue]  (2.22, 2.2) circle [radius = 0.04];
 \begin{scope}[white]
        \draw[fill=red!30, semitransparent] (7.5, 0) -- (7.5,2.5) -- (8,4) --  (8,1.5) --cycle;
        \draw[fill=red!30, semitransparent] (9, 0) -- (9,2.5) -- (9.5,4) --  (9.5,1.5) --cycle;
        \draw[fill=red!30, semitransparent] (10.5, 0) -- (10.5,2.5) -- (11,4) --  (11,1.5) --cycle;
\end{scope}
\draw[->, line width = 0.6, blue]  (7,0.5) to (12,0.5);
\draw[->, line width = 0.6, red]  (7.5,0) to (7.5,4);
\draw[->, line width = 0.2, red]  (7.4,0.2) to (8.2,2.6);
\node at  (12, 0.3) { \color{blue} \tiny$q$};
\node at  (7.3, 3.75) { \color{red} \tiny $w^1$};
\node at  (8.5, 2.5) { \color{red}  \tiny $ w^2$};
\draw [line width = 0.8, red, dashed](7.6,0.8)  to  [out=40, in=250] (7.9, 1.9)   ; 
\draw [line width = 0.8, red, dashed](7.9,2.2)  to  [out=87, in=270] (7.95, 3.8)   ; 
\draw[fill, red]  (7.9, 1.9) circle [radius = 0.04];
\draw[fill, red]  (7.9, 2.2) circle [radius = 0.04];
\draw [line width = 0.7, blue](7.6,0.8)  to  [out=-5, in=230] (9.2, 1.9)   ; 
\draw [line width = 0.7, blue](9.22,2.2)  to  [out=60, in=240] (11, 3.8)   ; 
\draw[fill, blue]  (9.2, 1.9) circle [radius = 0.04];
\draw[fill, blue]  (9.22, 2.2) circle [radius = 0.04];
\draw [line width = 1, red, densely dotted](7.9,1.9)  to   [out=80, in=270] (7.9, 2.2)   ; 
\draw [line width = 1, blue, densely dotted](9.2,1.9)  to [out=80, in=270] (9.22, 2.2)   ; 
 \end{tikzpicture}
 }
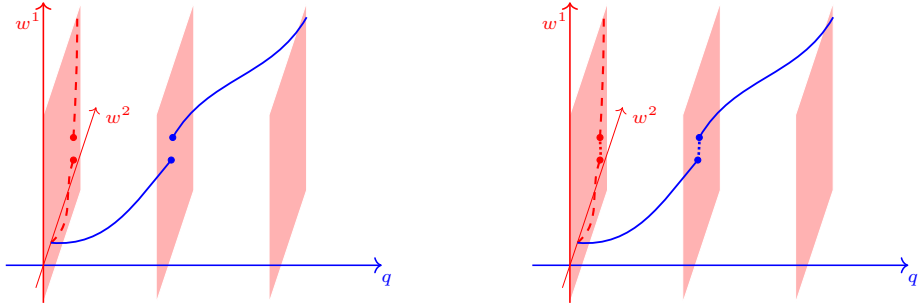
 \centerline{\tiny \bf \hskip - 0.3 cm Fig.1  $\cQ \times \cK$-projection of the graph of a solution \hskip 0.5cm  Fig.2 $\cQ \times \cK$-projection of a completed  graph }
Note that 
a   completed graph  can be also described as the (oriented) trace of a graph completion in the sense of Bressan and  Rampazzo 
(see  \cite{BR, GSZ1}).\par
\medskip
\subsection{Reachable sets and $\cM$-attainable sets}
 Given a state $q_o \in \cQ$ and   a positive time $T \in (0, + \infty)$,
 the {\it reachable set   in time  exactly equal to $T$ and by means of  piecewise $\cC^\o$  solutions}  is the subset of $\cQ$ defined by 
\begin{multline} \Reach^{\cC^\o}_{ T}(q_o) {:=}  \bigg\{ q \in \cQ\,:\, q\ \text{is the final point of a  piecewise $\cC^\o$ curve $q(t)$}\\
 \text{ that starts at $q_o$ and is the projection on} \ \cQ\\
  \text{of a piecewise $\cC^\o$ solution $(q(t), w(t))$ to \eqref{controlsystem} with}\  t \in [0, T]\   \bigg\}\ .
 \end{multline}
 The {\it $\cC^\o$-reachable set in time $T$} is  the set  
 $$\Reach^{\cC^\o}_{\leq T}(q_o) {=} \{q_o\} \cup \bigcup_{ T' \in (0, T]} \Reach^{\cC^\o}_{ T'}(q_o)\ .$$
 The control   system \eqref{controlsystem}  is said to have  the   {\it hyper-accessibility} (resp.  {\it accessibility in $\cC^\o$-sense})  {\it property} if  for any $q_o \in \cQ$ and $T> 0$, the set  $\Reach^{\cC^\o}_{ T}(q_o)$ (resp.  $\Reach^{\cC^\o}_{\leq T}(q_o)$) is open (resp. has non-empty interior).  It  has  the  {\it  small-time local controllability property  at $q_o \in \cQ$ in the $\cC^\o$ sense} if  there is a  $T> 0$ such that $ \Reach^{\cC^\o}_{ \leq T}(q_o)$  contains a neighbourhood of $q_o$. 
\par
Note that,  if the system has   the hyper-accessibility property, then it   is also accessible in the $\cC^\o$-sense. Moreover, {\it  for any point  $q_o$  with   the homing property} (i.e.   with the property that  there is   at least one  piecewise  $\cC^\o$-solution $(q(t), w(t))$ with   $q(0) = q_o = q(T)$,)   {\it the system has also  the  small-time local controllability property at $q_o$}.  The  simplest examples of points with  the homing property are  the  {\it  states  of equilibrium} (or {\it  state stable points} according to the terminology  of \cite{GSZ2}), i.e.   the   states $q_o$ for which there is  at least one choice for a control curve $w(t)$,  such that the corresponding  solution $t \mapsto (q(t), w(t))$ has $q(t) \equiv q_o$. 
\par
\smallskip
  The reachable sets of a control system are tightly related with the following subsets of  the extended space-time $\cM = \bR \times \cQ \times \cK$. Given  $x_o =  (t_o, q_o, w_o) \in \cM$,  the  {\it $\cM$-attainable set of $x_o$}  is the set 
 \begin{multline} \Attain_{x_o} {=} \bigg\{ y \in \cM\,:\, y\ \text{is the final endpoint of a  completed graph} \\
 \text{  starting from}\ x_o \ \text{and    corresponding to the graph of a piecewise}\ \cC^\o\ \text{solution}\ \bigg\}\ .
 \end{multline}
Let $\pi^\bR: \cM \to \bR$, $ \pi^\cQ: \cM \to \cQ$ be the standard 
projections of $\cM$ onto  its factors $\bR$ and $\cQ$, respectively. A set of reachable points $\Reach^{\cC^\o}_T(q_o)$, $T > 0$,   is  related with corresponding  $\cM$-attainable sets by 
\beq
\Reach^{\cC^\o}_T(q_o) = \bigcup_{\smallmatrix  w_o  \in \cK\\
 \endsmallmatrix}   \pi^\cQ \left(   \Attain_{(0, q_o, w_o)} \cap (\pi^{\bR})^{-1}( T)  \right) \ .
 \eeq
This relation yields  the following sufficient condition for hyper-accessibility:   {\it if for any $q_o \in \cQ$,  $T > 0$, there exists at least one  $w_o \in \cK$ such that  
$\pi^\cQ\bigg( \Attain_{(0, q_o, w_o)} \cap (\pi^\bR)^{-1}(T)\bigg) $
 is open in $\cQ$, then the system    is   hyper-accessible} and, in particular,  it  is small-time locally controllable at  its  points  with the homing property. \par
  \medskip
\section{Secondary distributions   and    $\cM$-attainable leaflets}
\label{sect3}
\subsection{The secondary distribution $(V^{II}, \cD^{II})$  of a control system} The  {\it canonical vector fields}  of the system \eqref{controlsystem} are  the vector fields  on $\cM = \bR \times \cQ \times \cK$  defined  by 
\beq \label{the-1} 
\bT^o    \= \frac{\p}{\p t} + f^i(t, q, w)\frac{\p}{\p q^i}\ ,\qquad
 W^o_{1}\= \frac{\p}{\p w^1}\ ,\qquad \ldots  \ ,\qquad  W^o_{m}\= \frac{\p}{\p w^m}\ .
 \eeq
%
%
For any system of this kind 
the {\it secondary distribution}    is the  pair  $(V^{II}, \cD^{II})$  defined by:
\begin{itemize}[leftmargin = 20pt]
\item[--] $V^{II}$ is the set of all (locally defined) real analytic  vector fields $X$, which are pointwise  finite linear combinations  of the   vector fields $W^o_\a$
and of  iterated Lie brackets 
\beq \label{genfor} W^{(k)}_\b \=  \underset{\text{$k$-times}}{\underbrace{[\bT^o, [\bT^o, \ldots, [ \bT^o, }}W^o_\b]\ldots]]\ , \qquad 1 \leq \b \leq m\ ,\ k = 1, 2, \ldots \ ;\eeq
\item[--] $\cD^{II}$ is the collection  of the spaces $\cD^{II}_x \subset T_x \cM$, $x \in \cM$,  spanned by the values at the point $x$ of the vector fields  in  $V^{II}$.
\end{itemize}  
We recall  that $(V^{II}, \cD^{II})$ is a  real analytic {\it generalised distribution}  in the sense of \cite[Def. 2.1]{GSZ2} (such a definition is later   recalled in  \S \ref{sect41}). When all spaces  $\cD^{II}_x$   have  an identical     dimension $p$, then $\cD^{II}$ is a regular distribution of rank $p$  in the classical sense and $V^{II}$ coincides with the set of all local vector fields  taking values in $\cD^{II}$. But there are cases in which the dimensions of the spaces $\cD^{II}_x$ are not constant and the collection of vector fields $V^{II}$ cannot be  recovered in a canonical way from $\cD^{II}$ as in the regular setting (see  \S \ref{sect41} for more detail).\par
The following examples give some  illustrations of the notion of the  secondary distributions. Others  examples are presented and  carefully discussed in \S \ref{applications}.
\begin{example} \label{example21}
 Consider a   linear control system,  i.e.   a system of the form
\beq \label{KAL0} \dot q^i = A^i_j q^j + B^i_\a w^\a\ ,\qquad A = (A^i_j) \in \bR_{n \times n}\ ,\ \ B = (B^i_	\a) \in \bR_{n \times m}\ . \eeq
In this case  the associated  canonical vector fields are 
\beq \label{KAL1} \bT^o  \= \frac{\p}{\p t}  + (A^i_j q^j + B^i_\a w^\a)\frac{\p}{\p q^i}\ ,\qquad W^o_\a \=  \frac{\p}{\p w^\a} \ , \ \ 1 \leq \a \leq m\ .\eeq
The  set  $V^{II}$ is indeed given by the  local real analytic  vector fields that  are pointwise linear combinations  of the vector fields 
 \begin{align}
 \nonumber &   W^{(0)}_1\=  W^o_1 \ ,\  W^{(0)}_2\=  W^o_2 \ , \ldots\ , && W^{(0)}_m \=  W^o_m \ ,\\
 \label{voila} & W^{(1)}_1 \=  [\bT^o, W^o_1] = - B^i_1\frac{\p}{\p q^i}\ , \ldots\ , &&W^{(1)}_m \=  [\bT^o, W^o_m] = - B^i_m\frac{\p}{\p q^i} \ ,\\
 \nonumber & W^{(2)}_1 \= [\bT^o, [\bT^o, W_1^o]] =   A^i_j B^j_1\frac{\p}{\p q^i}\ ,  \ldots\ , &&   W^{(2)}_m \=  [\bT^o, [\bT^o, W^o_m] ]=   A^i_j B^j_m\frac{\p}{\p q^i}\ ,\\
\nonumber &  & \text{etc.} &&  
 \end{align}
 The family $\cD^{II}$ is the collection of the  subspaces  $\cD^{II}_x \subset T_x\cM$ which are spanned by the values  of these vector fields at the points  of  $\cM$.   All spaces $\cD^{II}_x$ have  the same dimension and   $(V^{II}, \cD^{II})$ is a regular distribution. Note  that the vector fields $W_\a^{(0)}$,  $W_\a^{(1)}$,  $W_\a^{(2)}$, $ \ldots$, are in addition commuting generators for the regular   distribution $\cD^{II}$. This means that $\cD^{II}$ is also involutive.
 \end{example}
 \begin{example}  \label{example32} Consider now the control system with $\cQ = \bR$,  $\cK \subset \bR$ and  only one  control   equation  $\dot q = w^2$. In this case,  $\cM = \bR \times \cQ \times \cK \subset  \bR^3$,  the  canonical vector fields are  
\beq \label{KAL1*} \bT^o  \= \frac{\p}{\p t}  +  w^2\frac{\p}{\p q}\ ,\qquad W^o \=  \frac{\p}{\p w}\eeq
and the  set $V^{II}$ consists of all local real analytic   vector fields that are  pointwise linear combinations   of  the vector fields
\beq W^{(0)} = W^o = \frac{\p}{\p w}\ ,\qquad   W^{(1)} =  [\bT^o, W^{(0)}] = - 2  w\frac{\p}{\p q}\ .\eeq
Indeed, since $[\bT^o, W^{(1)}] = 0$,  there are no other non trivial vector fields  of the form  \eqref{genfor}. In contrast with the previous example,  
the   family of vector spaces  $\cD^{II}$ is not a regular distribution, since   the dimension of   $ \cD^{II}_{x = (t, q, w)} = \langle  X_{x= (t,q,w)}, X \in V^{II}\rangle$ is   $1$  when  $w = 0$ and $2$ otherwise (on this secondary distribution, see  also \S \ref{elementary}).
\end{example}
\par
\medskip
\subsection{Surrogate  fields and surrogate leaflets  centred at $\cM$-attainable points} \label{section32}
Let   $x_o = (t_o, q_o, w_o) \in \cM$, $y_o = (\wh t_o, \wh q_o, \wh w_o)  \in \Attain_{x_o}$  and   $\h = \h_1 \ast \h_2 \ast  \ldots\ast \h_{2r+1}$ a  completed graph  in  $\cM = \bR \times \cQ \times  \cK$ that starts from $x_o$ and ends at $y_o$.  
 We may  assume that  the piecewise regular curve $\h$ is decomposed into a sufficiently large number of small regular arcs, so that  (a) the last arc  is entirely  contained in a sufficiently small prescribed neighbourhood $\cU \subset \cM$ of $y_o$ and (b)   such a  last arc is parametrisable by a real analytic map  $s \mapsto \h_{2r +1}(s)$ with  $s$ running in a  small interval $[0, \ve]$. 
A remark on notation:  As in \cite{GSZ1}, in all what follows write  $\Phi^X_s$  for  the $s$-parameterised family of local diffeomorphisms which constitute the  flow of  a vector field $X$.\par
\smallskip
Applying    \cite[Rem.3.6]{GSZ2}  to the curves $\h_{2r +1}(s)$ and  $\wh \h(s) = \Phi^{\bT_o}_{s-\ve}(y_o)$ (which are both ending at  $y_o = \h_{2r +1}(\ve) = \wh \h(\ve)$), we conclude that    there is no loss of generality if we assume that   $\h_{2r +1}(s)$ has the form 
$$\h_{2r + 1}(s) = \Phi^{\bT}_{s}(x_o')\ ,\qquad s \in [0, T] $$
for some point $x_o'$ sufficiently close to $y_o$, and a pushed-forward vector field $\bT \= \f_*(\bT^o)$  by a real analytic  local diffeomorphism $\f$ of the form $\f(t, q, w) = (t, q, \wh \f(t,q,w))$.  Due to  the particular form of the diffeomorphism $\f$, the difference between the vector fields $\bT^o$ and $\bT = \f_*(\bT^o)$  has the form
$$\d \bT = \bT - \bT^o =  \bT^o(\wh \f^\a)  \frac{\p}{\p w^\a}$$
and has therefore  trivial components along the vector fields $\frac{\p}{\p t}$ and $\frac{\p}{\p q^i}$. \par
\smallskip
The next Theorem  \ref{thm34}  can be taken as  the foundation of  the two criterions which we are going to present in the next section. For stating it, we  need to introduce some new objects (\cite[Def. 5.8 and Def. 7.4]{GSZ2}).  The vector fields  on the neighbourhood  $\cU$ of $y_o$ of the form 
\beq W^{\t}_\a = \Phi^{\bT}_{\t*}(W_\a) = W_\a + \sum_{k = 1}^\infty \frac{(-1)^k}{k!} \underset{\text{$k$-times}} {\underbrace{[\bT, [\bT, \ldots, [\bT, W_\a]\ldots]]}}\t^k\ , \qquad \t \in (0, T)
\eeq
are called {\it elementary $\bT$-surrogate fields  of $\bT$-depths $\t$}. Any  linear combination  $X = \lambda^\a W_\a^{\t} $  of elementary $\bT$-surrogate fields $W_\a^{\t} $, all with the same $\bT$-depth,  is called {\it $\bT$-surrogate field} of  $\bT$-depth $\t$.\par
\smallskip
 \begin{definition}  \label{def33}
 Given a vector field  $\bT = \f_*(\bT^o)$ as above, 
 a {\it  $\bT$-surrogate map of rank $M$,  centred at $y_o$ and  with $\bT$-depths in $(0, T)$}  is a smooth  map $F: \cV \subset \bR^M \to  \cM$  of class $\cC^\o$ from a neighbourhood $\cV = (-\ve, \ve)^M$ of the origin $0 \in \bR^M$  
satisfying  the following: 
\begin{itemize}[leftmargin = 20pt]
\item[(1)] $F(0) = y_o$ and  $F(\cV)$ is an embedded $M$-dimensional submanifold of $\cM$; 
\item[(2)] there  is 
\begin{itemize}
\item a set of  $\bT$-surrogate  fields $X_1, \ldots, X_{\bf m}$, $M \leq {\bf m}$, with strictly decreasing $\bT$-depths $\t_\ell$  in $(0, T)$
$$T > \t_1 > \t_2 > \ldots > \t_{\bf m}> 0\ ,$$
 \item a set of smooth  functions  $\s_\ell: (- \ve, \ve) \subset \bR \to \bR$,  
 \item an ${\bf m}$-tuple $(i_1, \ldots, i_{\bf m})$ of integers, each one taken  in the set $\{1, \ldots, M\}$ 
 \end{itemize}
  such that for any $(s^i) \in \cV$  the corresponding point  $F(s^1, \ldots, s^M) \in \cM$ is equal to
\beq \label{cond1} F(s^1, \ldots s^M) = \Phi^{X_{\bf m}}_{\s_{\bf m}(s^{i_{\bf m}})} \circ \ldots\circ  \Phi^{X_1}_{\s_1(s^{i_1})}(y_o) \ .\eeq
\end{itemize} 
The images of  $\bT$-surrogate maps are called  {\it $\bT$-surrogate} (or just {\it surrogate}) {\it  leaflets}. 
\end{definition}
The compositions of flows, which  define surrogate leaflets,  are  those which allow to pass directly  from a given point $y_o \in \Attain_{x_o}$ to  other points  in $\Attain_{x_o}$, as  we  mentioned in the Introduction. For a detailed description of the geometric properties of the $\bT$-surrogate fields and of the surrogate leaflets, we refer  to \cite{GSZ2}. For the purposes of this paper, we just need   the  following   consequence of the  discussion  in  \cite[\S 7.2]{GSZ2}.
 \begin{theo} \label{thm34} Let  $x_o \in \cM$ and $y_o \in \Attain_{x_o}$. Then there is a connected neighbourhood $\cU \subset \cM$ of $y_o$, a vector field $\bT = \f_*(\bT^o)$,  with  $\f$  of the form $\f(t, q, w) = (t, q, \wh \f(t,q,w))$, and  a $T> 0$, such that  all  $\bT$-surrogate leaflets, centred  at $y_o$ and   with   $\bT$-depths in $(0, T)$,    are entirely included in  $ \Attain_{x_o}$. 
\end{theo}
\par 
\medskip
\section{The basic relations between   $\cM$-attainable sets, surrogate leaflets and  $\cD^{II}$-path connected components of $\cM$}
\label{sect4}
As we mentioned in the Introduction, the surrogate leaflets   are important not only because they  consist of $\cM$-attainable points  from the same point $x_o$ of their center,  but also because they  are also tightly related with the orbits of the pseudo-group $\cG^{II}$ of local diffeomorphisms,  which is generated by the flows of the vector fields in the secondary distribution   $(V^{II}, \cD^{II})$. The dimensions of such orbits  are in many cases  computable with the help of Chow-Rashevski\u\i-Sussmann's Theorem and Frobenius Theorem. In this section we  brush up with  this results  and, immediately after, we discuss in greater detail their relations with the  surrogate leaflets.\par
\smallskip
\subsection{Points  joined   by  paths  tangent to a distribution: Chow-Rashevski\u\i-Sussmann Theorem} \label{sect41}
In \cite[Def.2.1]{GSZ2},  we introduced the  notion of   ``generalised distribution',  which is essentially just a variant  of many  other classical definitions of  non-regular distributions on manifolds -- see e.g.  \cite{Na, SJ, Su}. 
 A real analytic {\it generalised distribution}  on a manifold $\cM$  is  a pair $(\wc V, \wc \cD)$, consisting of:
\begin{itemize}[leftmargin = 15 pt]
\item  a set   $\wc V$  of vector fields, which is  locally generated over the ring  of  the real analytic functions  by  finite sets of real analytic  vector fields $(X_1, \ldots, X_p)$ and which is constrained by the following requirement:  if $(X_1, \ldots, X_p)$ and $(X'_1, \ldots, X_{p'}')$   are two tuples of local generators on the same open set $\cU \subset \cM$, one tuple  can be expressed in terms of the other by a  matrix with real analytic functions as  entries; 
\item  the collection $\wc \cD$ of all spaces $\wc \cD_x \subset T_x \cM$  which are spanned by the values at $x$ of the vector fields in the set $\wc V$.
\end{itemize}
If the  spaces $\wc \cD_x$ have equal  dimension for any $x \in \cM$,  the family  $\wc \cD$   is a regular  distribution and the set of vector fields  $\wc V$  coincides with the set of all local $\cC^\o$ vector fields  taking values in $\wc \cD$.  In this case, we say that  the pair $(\wc V, \wc \cD)$  is  {\it regular}.\par
\smallskip
Given a  $\cC^\o$ generalised distribution  $(\wc V, \wc \cD)$  on  $\cM$, a piecewise $\cC^\o$ regular curve $\g = \g_1 \ast \g_2 \ast \ldots \ast \g_r$  is  said to be    {\it tangent to $(\wc V, \wc \cD)$}  (or   {\it $\wc \cD$-path},  for short)  if all tangent vectors of the regular arcs $\g_i$ are  in the spaces $\wc \cD_x$, $x \in \cM$.   The equivalence classes in $\cM$ for  the  relation 
$$y \sim y'\qquad \Longleftrightarrow \qquad y \ \text{and} \ y' \ \text{are joined by a $\wc \cD$-path}$$
are called {\it $\wc \cD$-path components of $\cM$}. \par
\smallskip
Note that {\it the $\wc \cD$-path component of a point $x_o$  coincides with   the  orbit $\cG^{\wc V}(x_o)$ of the  pseudogroup $\cG^{\wc V}$ of local diffeomorphisms, which is determined by  all   finite compositions of local   diffeomorphisms of  the flows of the  vector fields in  the set $\wc V$}.  This  can be checked as follows.  By definition, any point $x$ which can be reached from a given $x_o \in \cM$ applying  a  finite composition of local diffeomorphisms of  flows of vector fields in $\wc V$, is equivalent to $x_o$ by  the above equivalence relation. This  means that the orbit   $\cG^{\wc V}(x_o)$  is  contained in the  $\wc \cD$-path component of the  point $x_o$.  On the other hand, by Sussman's Theorem  \cite{Su}, any  orbit $\cG^{\wc V}(x_o)$ is an immersed submanifold  of $\cM$ with each  tangent space $T_x \cG^{\wc V}(x_o) \subset T_x \cM$,   which includes the corresponding  space $\wc \cD_x$. In particular,  it can be covered by  a  countable  union of  embedded submanifolds $\cV \subset \cM$,  whose   tangent spaces include all  spaces $\wc \cD_x$, $x \in \cV$. This implies that any $\wc \cD$-path  starting from  a point $x_o$  is  entirely included in the immersed submanifold   $\cG^{\wc V}(x_o)$, meaning that  the $\wc \cD$-path component  of $x_o$ is a subset of $\cG^{\wc V}(x_o)$, in fact equal to it. \par
\smallskip
The above mentioned Sussmann's Theorem, combined with the Noetherian property of the  rings of  real analytic  functions  (\cite[Thm. 3.8]{Ma}),  has several useful consequences in  case of  the real analytic 
generalised distributions (see  \cite[Thm. 5.16 and Cor. 5.17]{AS}).  We briefly  recall  them in the next theorem and in the   subsequent discussion (see also \cite[\S 7.1]{GSZ1}. In what follows,  given a finite set   $\{Y_1, \ldots, Y_m\}$ of local vector fields defined  on a common  open subset of $\cM$, for any $r \geq 2$ we denote by $Y_{(i_1, \ldots, i_r)}$ the iterated Lie bracket
$$Y_{(i_1, \ldots, i_r)} \= [Y_{i_1}, [Y_{i_2}, [\ldots [Y_{i_{r-1}}, Y_{i_r}]\ldots]]]\ ,\qquad 1 \leq i_\ell \leq m\ .$$
For $r = 1$, we  set $Y_{(i_1)} \= Y_{i_1}$.  The integer $r$ is called  the {\it depth} of the iterated Lie bracket.
\begin{theo} \label{Rash} Let   $(\wc V, \wc \cD = \cD^{\wc V})$ be a real analytic generalised distribution and denote by  $(\wc V^{(\text{\rm Lie})}, \wc \cD^{(\text{\rm Lie})})$ the pair given by the family $\wc V^{(\text{\rm Lie})}$  of  all real analytic vector  fields that are finite combinations  of  vectors $Y_{(i_1, \ldots, i_r)}$, with vector fields $Y_j \in \wc V$,  and the associated family of spaces  $\wc \cD^{(\text{\rm Lie})}_y \subset T_y \cM$,  spanned by the vector fields in  $\wc V^{(\text{\rm Lie})}$. Then:
\begin{itemize}[leftmargin = 20pt]
\item[(i)] The pair  $(\wc V^{(\text{\rm Lie})}, \wc \cD^{(\text{\rm Lie})})$ is an involutive  real analytic  generalised distribution; 
\item[(ii)] The $\wc \cD$-path connected component of a  point $x_o$ is  the  maximal integral leaf $\cS$  through $x_o$ of  the  bracket generated distribution  $(\wc V^{(\text{\rm Lie})}, \wc \cD^{(\text{\rm Lie})})$ and 
 any such maximal  integral leaf  $\cS$ is an immersed submanifold of $\wc \cM$;  
 \item[(iii)] Any $\wc \cD$-leaflet (= the  image of a map of the form \eqref{cond1} for vector fields in $\wc \cD$) is included in  a  unique maximal integral leaf $\cS$ of  $(\wc V^{(\text{\rm Lie})}, \wc \cD^{(\text{\rm Lie})})$ and for  any point $x_o$ of  such integral leaf, there exists a neighbourhood $\cU$  such that   $\cS \cap \cU$   is a $\wc \cD$-leaflet of maximal dimension.
\end{itemize}
\end{theo}
\par
 We recall that a maximal integral leaf $\cS^{(x_o)}$  of  $(\wc V^{(\text{Lie})}, \wc \cD^{(\text{Lie})})$ is  an immersed submanifold, i.e. it is the image $\imath^{(x_o)}(S)$ of  a smooth immersion $\imath^{(x_o)}: S \longrightarrow \cM_j$  of a manifold $S$ of dimension $\dim S = \dim \cD^{(\text{Lie})}_{x_o}$. In  case the immersion is not an embedding, the topology on  $\cS^{(x_o)}$,  given by  the images under $\imath^{(x_o)}$ of the open sets of  $S$, does not coincide with  the  topology of $S^{(x_o)}$  as a subset of $\cM$. In order to distinguish these two  topologies,   we  call  the first    the  {\it  intrinsic topology}  of  $\cS^{(x_o)}$
\par
Now, given a real analytic generalised distribution    $(\wc V, \wc \cD = \cD^{\wc V})$ and denoting by  $(\wc V^{(\text{\rm Lie})}, \wc \cD^{(\text{\rm Lie})})$ the  corresponding  bracket generated  distribution,  a  {\it  decomposition of  an open subset $\cU \subset \cM$ into  $\wc \cD$-strata of  maximal $\wc \cD$-depth $\mu$} is a finite family of disjoint subsets $\cU_0$, $\cU_1$, \ldots, $\cU_p$   such that $\cU = \cU_0 \cup \cU_1 \cup \ldots \cup \cU_p$  and  
\begin{itemize}[leftmargin = 20pt]
\item[(i)] for each $\wc \cD$-stratum $\cU_j$,    the spaces   $\wc \cD^{(\text{Lie})}|_{y} \subset T_y\cM$, $y \in \cU_j$,   have all the same dimension  and all  maximal integral leaves in $\cU$ of   $(\wc V^{(\text{\rm Lie})}, \wc \cD^{(\text{\rm Lie})})$ passing through the points of $\cU_j$ are  entirely  included  in $\cU_j$; 
\item[(ii)] there is an integer $\mu \geq 1$ (called {\it maximal $\wc \cD$-depth}) and a set of integers $1 \leq \mu_j \leq \mu$,  one per each stratum $\cU_j$, such that  for each $\cU_j$ the spaces  $\wc \cD^{(\text{Lie})}|_{y}$, $y \in \cU_j$, are  generated by the values at $y$ of a finite number of  iterated Lie brackets  $Y_{(i_1, \ldots, i_r)}$,  with  $Y_{j_\ell}\in \wc \cD|_{\cU_j}$ and depth  $r \leq \mu_j$.   
\end{itemize}
Note  that, by  \cite[Thm. 3.8]{Ma}, for any   $x_o \in \cM$,  there is  a neighbourhood $\cU$ such that  $\wc \cD^{(\text{Lie})}|_{\cU}$ is spanned by  a finite number of  iterated Lie brackets  $Y_{(i_1, \ldots, i_r)}$. Then a  decomposition  of $\cU$ into  $\wc \cD$-strata can be  determined as follows:  Let  $\cU_0 \subset \cU$ be   the maximal set of points  $y \in \cU$   for which $\dim \wc \cD^{(\text{Lie})}|_{y}$ is maximal;  Then,  set   $\cU_1$  to  be   the maximal  set of  points  $y \in \cU \setminus \cU_0$  for which     $\dim \wc \cD^{(\text{Lie})}|_{y}$ has the second  maximal value,    and so on.  By construction and Theorem \ref{Rash}, each maximal integral leaf of $\wc \cD^{(\text{Lie})}|_{y}$ through a point of $y \in \cU_j$ is necessarily entirely  included in  $\cU_j$ and,  being  $\wc \cD|_{\cU}$  spanned by  finite number of iterated  Lie brackets,  integers  $1\leq \mu_j \leq \mu = \max \mu_j$  for which (ii) holds can be directly determined. \par

\begin{example} \label{example44} We have seen that  the secondary distribution $(V^{II}, \cD^{II})$  of  Example \ref{example21} is regular and involutive. Thus, according to the above recipe for  determining $\cD^{II}$strata,   the whole manifold  $\cM$ admits a  decomposition into $\cD^{II}$-strata, with just one  stratum, namely $\cU_0 \= \cM$ itself.  The  bracket generated distribution $\cE = \cE^{(\cU_0)} $ is $\cE =  \cD^{II}$. The rank of  $\cE$ is   the maximal number of pointwise linearly independent vector fields   in the set 
$$\bigg\{ W^{(0)}_1 ,\ \ldots \ , W^{(0)}_m, W^{(1)}_1 ,\ \ldots\ ,\ W^{(1)}_m\ ,\  W^{(2)}_1 \ , \ldots,  W^{(2)}_m \ , \ldots\ \bigg\}\ .$$
This number  is also   the maximal number of linearly independent vectors in the following set of $n$-vectors 
(here, for   $v = (v^j) \in \bR^n$,  we denote by   $A {\cdot} v$ the standard multiplication between  the squared matrix $A = (A_i^j)$ and the column $v = (v^j)$):
 \begin{multline} \label{tre}B_1, B_2 \ldots, B_m, \ \ A{\cdot}B_1, \ldots, A{\cdot} B_m,\ \  (A{\cdot} A){\cdot}  B_1,  (A{\cdot} A){\cdot}  B_1 \ldots, \ (A {\cdot} A){\cdot} B_m, \ldots\\
\ldots, \ \ \left(A{\cdot} \ldots {\cdot}A\right){\cdot}  B_1,  \left(A {\cdot}\ldots {\cdot}A\right){\cdot}  B_2, \ldots, \left(A{\cdot} \ldots {\cdot} A\right) {\cdot} B_m,\ \ \ldots
 \end{multline}
The maximal integral leaves of the distribution $\cE$ (= the $\cD^{II}$-path components of $\cM$) are the intersections of  affine subspaces of $\bR^{1 + n + m}$ of dimension $p = \rank \cE$ with  $\cM = \bR^{1 + n} \times \cK \subset \bR^{1 + n + m}$.  Each of them  is  contained in a hypersurface  of the form $\{t = \text{const.}\}$
\end{example}
\par
\begin{example} \label{example45} The secondary distribution $(V^{II}, \cD^{II})$ of  the Example \ref{example32}  is not regular. On the other hand,  for any $x\in \cM$,  the  space $E^{\text{Lie}}_{x} \subset  T_x \cM$ which is spanned by the values in $x$ of the vector fields  in $V^{II}$ and of their   iterated Lie brackets  is    
 \begin{multline*} E^{\text{Lie}}_{x= (t, q, w)}  = \left\langle W^{(0)}_x = \frac{\p}{\p w}\bigg|_x\ ,\ W^{(1)}_x =  -2  w \frac{\p}{\p q}\bigg|_x\ , [W^{(0)}, W^{(1)}]_x = - \frac{\p}{\p q}\bigg|_x \right\rangle = \\
 =  \left\langle \frac{\p}{\p w}\bigg|_x\ ,\ \frac{\p}{\p q}\bigg|_x \right\rangle\ .
 \end{multline*}
 Hence,  also in this case, the whole $\cM$ admits a decomposition into $\cD^{II}$-strata,  consisting of just one stratum $\cU_0 \= \cM$.  And one can immediately  check that  the  bracket generated  distribution $\cE = \cE^{(\U_0)}$  is spanned at all points by  the commuting vector fields $\frac{\p}{\p q}$ and $\frac{\p}{\p w}$. Its maximal integral leaves   (i.e.  the $\cD^{II}$-path  components of $\cM$)  are  the intersections of the hyperplanes $\{t = \text{cost.}\}$ of $\bR^3$ with  $\cM = \bR^2 \times \cK  \subset \bR^3$.
 \end{example}
\par\smallskip
 \subsection{$\cM$-attainable sets and $\cD^{II}$-path  components}
 We are now ready to state the main results of \cite{GSZ2} concerning  the relations between the   $\cG^{II}$-orbits (= $\cD^{II}$-path  components) and the $\cM$-attainable sets. We first need the following proposition, which  is  a direct consequence of the definition of surrogate leaflet and of Theorem \ref{Rash}. \par
 \begin{prop} \label{theorem4.6}  Let $x_o \in \cM$,   $y_o \in \Attain_{x_o}$ and $\cU$ a neighbourhood of $y_o$ admitting a decomposition into $\cD^{II}$-strata $\cU = \cU_0 \cup \ldots \cup \cU_r$. Denote by  $\cU_{j_o}$ and $\cS^{(y_o)} \subset \cU_{j_o}$ the $\cD^{II}$-stratum  and the $\cD^{II}$-path  component of $\cU$, respectively,  which  contain  $y_o$.  
 For any sufficiently small  neighbourhood $\cU' \subset \cU$ of $y_o$, any sufficiently small $T> 0$, and any vector field $\bT = \f_*(\bT^o)$ on $\cU$,  with  $\f$ local diffeomorphism of the form $\f(t, q, w) = (t, q, \wh \f(t,q,w))$,  each  $\bT$-surrogate leaflet $\cL \subset \cU'$, centred at $y_o$ and  $\bT$-depths in $(0, T)$,  is  included in  $\cS^{(y_o)} \cap \cU'$.
\end{prop}
Combining Proposition \ref{theorem4.6}, Theorem \ref{thm34} and  the Implicit Function Theorem, the following  corollary follows (see also \cite[\S 7.2]{GSZ2}). \par
\begin{cor} \label{thecor}  Assume that for any  $y_o = (t_o, q_o, w_o) \in \cM$ and  any 
local vector field $\bT = \f_*(\bT^o)$ with $\f$
local diffeomorphism $\f: \cU \subset \cM \to \cM$ around $y_o$ of the form $\f(t, q, w) = (t, q, \wh \f(t,q,w))$  the following two conditions hold:  Denoting by  $\cS^{(y_o)}$  the $\cD^{II}$-path  component  of $y_o$, 
\begin{itemize}[leftmargin = 20pt]
\item[(1)] there exists a  $\bT$-surrogate leaflet $\cL$, which is centred at $y_o$ and is open  in $\cS^{(y_o)}$, 
\item[(2)]  the standard projection $\pi^\cQ: \cM = \bR\times \cQ \times \cK \to \cQ$ maps homeomorphically a neighbourhood of $y_o$  of the intrinsic topology of  $\cS^{(y_o)}$   onto   a neighbourhood   of $q_o$. 
\end{itemize}
 Then the  control system \eqref{controlsystem} has the hyper-accessibility property.
\end{cor} 
The points $y_o \in \cM$ for which  (1) holds are called {\it good points} (\cite[Def. 7.2]{GSZ2}). If all points  are good, we say that the control system is  {\it essentially good}. With this terminology, the   corollary  can be   stated as follows.   {\it Assume that   the control system is essentially  good. If the restriction of the   projection $\pi^\cQ$   to any   $\cD^{II}$-path component  $\cS^{(y_o)}$ of $\cM$ has maximal rank, then the control system has the hyper-accessibility property.}
\par
\smallskip
In most cases,  checking whether   all   restrictions   $\pi^\cQ|_{\cS^{(y_o)}}$  satisfy the maximal rank condition    reduces to   computing   the ranks of the   restricted  Jacobians  $J(\pi^\cQ)\big|_{\cE^{(\cU_j)}_x = (t, q, w)}$, i. e. the restrictions of the Jacobians   to the spaces of the  bracket generated distributions of the $\cD^{II}$-strata of a neighbourhood $\cU$ of $y_o$.  From a computational point of view,  this  is in principle  an affordable  task. 
What is not  trivial to check   is whether  the control system is essentially good or not.  In  the next subsection  \S \ref{thetwocrit}, we give two   efficient  criterions  for this purpose. \par
\medskip
\subsection{Two criterions for  goodness}\label{thetwocrit}
Our first criterion for  goodness is very simple:
\begin{theo}[{\cite[Thm. 7.9]{GSZ2}}]  \label{prop93}Let  $ y_o \in \cM$,  assume that  $\cU$ is a neighbourhood of $y_o$ admitting a decomposition into $\cD^{II}$-strata and denote by $\cU_{j_o} \subset \cU$   the $\cD^{II}$-stratum of $\cU$ containing $y_o$. If the restricted  distribution  $(V^{II}|_{ \cU_{j_o}}, \cD^{II}|_{ \cU_{j_o}})$  is regular and involutive near $y_o$, then   $y_o$  is a good point.
\end{theo}
 The points  in  Theorem \ref{prop93} are called {\it good points of the first kind}.\par
 \smallskip
 Our second criterion for goodness has  apparently a more involved statement. It is nonetheless equally  simple to be used, as it is  illustrated by the examples of  the next section. In order to state this criterion, we first  need to introduce the following definition (for more details, see  \cite[\S 5.1]{GSZ2}). \par
 \begin{definition}  \label{lemmone}   A   tuple $(W_A)$ of real analytic vector fields on an open set  $\cU\subset \cM$, which    generates   the secondary distribution $(V^{II}|_{\cU}, \cD^{II}|_{\cU})$,  is called {\it set of $\bT_o$-adapted generators}  if there is an integer $\nu \geq 0$ and  a set of   integers  $R_a\geq 1$ for any $0 \leq a \leq \nu$,  such that the family  of the indices  $\{A\}$  of the tuple  is in bijection with the set of indices 
 $$\bigg\{``\ell(a)j"\ ,\ 0 \leq \ell \leq a\ ,\ 1 \leq j \leq R_a\ ,\  0 \leq a \leq \nu\bigg\}\ ,$$
 so that the  vector fields $W_A$  can be ordered  according to the following table   \\
\!\!\resizebox{0.98\hsize}{!}{\vbox{
 \begin{align*} & W_{0(0)1}, \ldots, W_{0(0)R_{0}},  &&W_{0(1)1}, \ldots, W_{0(1)R_{1}},  &&   \ldots, && \ W_{0(\nu -1)1}, \ldots, W_{0(\nu-1)R_{\nu-1}}, &&W_{0(\nu)1}, \ldots, W_{0(\nu)R_{\nu}}, \\[10pt]
&  & &W_{1(1)1}, \ldots, W_{1(1)R_{1}},  && \ldots && W_{1(\nu-1)1}, \ldots, W_{1(\nu-1)R_{\nu -1}},&&W_{1(\nu)1}, \ldots, W_{1(\nu)R_{\nu}}, \\[10pt]
 & && &&    && \ddots && \vdots && \\[10pt]
   & && && && W_{\nu-1(\nu -1)1}, \ldots, W_{\nu-1(\nu-1)R_{\nu-1}}, && W_{\nu-1(\nu)1}, \ldots, W_{\nu-1(\nu)R_{\nu}}, \\[10pt]
   & && && && &&  W_{\nu(\nu)1}, \ldots, W_{\nu(\nu)R_{\nu}}\ , 
   \end{align*}}}
   \\
and such that  the following three conditions hold:
\begin{itemize}[leftmargin = 20pt]
\item[(1)] The vector fields
$$W_{0(0)1}\ , \ \ldots\ , \ W_{0(0)R_0}\ , W_{0(1)1}\ , \ \ldots\ , \ W_{0(1)R_1}\ ,\ \ldots\ ,\ W_{0(\nu)1}\ , \ \ldots\ , \ W_{0(\nu)R_\nu}$$ 
are  pointwise linearly independent generators for   the same distribution, which is spanned by the canonical vector fields $W_\a^o$  (in most cases, up to a reordering,  they are nothing but     the  $W^o_\a$);   in particular,  $R_0 + R_1 + \ldots + R_\nu =  m$; 
\item[(2)]  For  any   $1\leq a\leq \nu$, $0 \leq \ell \leq a$ and $1 \leq j \leq R_{a}$,  the corresponding vector field $W_{\ell(a), j}$  is obtained  from the $W_{0(a)j}$ by applying $\ell$-times the operator $X \mapsto [\bT^o, X]$, i.e. 
\begin{multline} \label{defW} W_{1(a)j} = [\bT^o, W_{0(a)j}]  \ ,\ \  W_{2(a)j} =  [\bT^o, [\bT^o, W_{0(a)j}]]\ ,\ \ \ldots \ \\[3 pt]
\ldots \ W_{\ell(a)j} =  \underset{\ell\text{\rm -times}}{\underbrace{[\bT^o, [\bT^o, [\bT^o,  \ldots [\bT^o, }}W_{0(a)j}] \ldots ]]]\ , \ \ \ldots \ \  \\
\ldots\ , \ \ W_{a(a)j} =  \underset{a\text{\rm -times}}{\underbrace{[\bT^o, [\bT^o, [ \bT^o, \ldots [\bT^o, }}W_{(0(a)j}] \ldots ]]] \ .
\end{multline}
 \item[(3)] For any   $0 \leq a\leq \nu$ and $1 \leq j \leq R_{a}$ the iterated Lie bracket
$$   \underset{(a+1)\text{\rm -times}}{\underbrace{[\bT^o, [\bT^o, [\bT^o,  \ldots [\bT^o, }}W_{0(a)j}] \ldots ]]] $$
is pointwise a linear combination of the  generators $W_{\ell'(b)j'}$ with   $b \leq a$.
\end{itemize}
 \end{definition}
By \cite[Lemma 5.4]{GSZ2},   for  any sufficiently small open set  $\cU$  there is a set of $\bT^o$-adapted generators on $\cU$. Explicit examples which illustrate this  notion  can be found  in  \S \ref{applications}. \par
\smallskip
We now recall that, according  to the terminology of \cite[\S 7.4]{GSZ2}, 
for any  canonical vector field  $W^o_\a  = \frac{\p}{\p w^\a}$,  the  {\it  secondary sub-distribution  generated by  $W^o_\a$} is the     generalised distribution  $(V^{II(W^o_\a)}, \cD^{II(W^o_\a)})$,  which has the  vector field $W^o_\a$ and all possible iterated Lie brackets  
$[\bT^o, [\bT^o, \ldots[\bT^o, W^o_\a]\ldots]]$ as generators.\par
\smallskip
We are now ready to state our second criterion.\par
\begin{theo}[{\cite[Thm. 7.13]{GSZ2}}] \label{criterione}    Let  $y_o \in \cM = \bR \times \cQ \times \cK$,  assume that  $\cU$ is a neighbourhood of $y_o$ admitting a decomposition into $\cD^{II}$-strata and denote by $\cU_{j_o} \subset \cU$   the $\cD^{II}$-stratum containing $y_o$.  Then $y_o$ is a good point if  the following three conditions are satisfied:  
\begin{itemize}[leftmargin = 15pt]
\item[(1)] the bracket generated distribution $\cE^{(\cU_{j_o})}$ has stratum depth $\mu_{j_o} = 2$;  
\item[(2)] denoting by ${\bf m}$ the rank of  $\cE^{(\cU_{j_o})}$, there exists  a neighbourhood $\cU' \subset \cU$ of $y_o$  and a set of  generators for $\cE^{(\cU_{j_o})}|_{\cU}$
$$\left(W_J\right)_{1 \leq J  \leq {\bf m}} \ ,$$ 
in which each vector field $W_J$ is one of the following two types: 
\begin{itemize}[leftmargin = 15pt]
\item[(a)] it  is an element $W_J = W_{A_J}$ of  a set $(W_A)$ of $\bT^o$-adapted generators, 
\item[(b)] it is a Lie bracket $W_J =  [W^{(W^o_\b)} _{B}, W_{B'}]$ given by 
\begin{itemize}[leftmargin = 10pt]
\item[--] a vector field $W^{(W^o_\b)} _{B}$ which  is an element of a set of  $\bT^o$-adapted generators of  a  secondary sub-distribution $(V^{II(W^o_\b)},\cD^{II(W^o_\b)})$,  for which   $y_o$ is a good point of the first kind; 
\item[--] a vector field  $W _{B'}$ which  is one of the $\bT^o$-adapted generators   for $(V^{II}, \cD^{II})$ of the  tuple  $(W_A)$; 
\end{itemize}
 \end{itemize}
 \item[(3)] At least one of the vector fields $W_J$ has the form (b).
 \end{itemize}
\end{theo}
 The points  for which   Theorem \ref{criterione} holds are called {\it good points of the second kind}. \par
 \smallskip
 As a concluding remark of this section, we would like to point out that, as it is mentioned in the Introduction of \cite{GSZ2}, the proof of this criterion  is structured in such a way  that it should not be hard to be  generalised to  settings, in which the  stratum depth  $\mu_{j_o} $ is  greater than or equal to $3$. This very plausible development is left to future investigations. \par
 \smallskip
\par
\medskip
\section{Applications} \label{applications}
\subsection{A new  proof of the Kalman criterion} 
\label{sectKal} 
We now give a new proof of the well known Kalman Theorem. The traditional  proof for  this classical result can be  found in most of the books on Control Theory --  see for instance  \cite{Bl, Co, Ju}. For our  proof, we   need the following preliminary result, which is of some interest of its own.
\par
\begin{theo} \label{theo51}  A linear control system  $ \dot q = A q  + B w$, with     $A = (A^i_j) \in \bR_{n \times n}$,  $B = (B^i_a) \in \bR_{n \times m}$ and  $\cK$ open subset of  $ \bR^m$,  has the hyper-accessibility property if and only if   the 
maximum $n_{\max}$ of the  ranks  
\begin{multline} \label{seq}
 n_\ell \= \dim \bigg\langle B_1, \ldots, B_m, \ \ A{\cdot}B_1, \ldots, A{\cdot}B_m, \ \ A {\cdot} A{\cdot}B_1, \ldots, A {\cdot} A{\cdot}B_m, \\
 A {\cdot} A {\cdot} A{\cdot}B_1, \ldots, A {\cdot} A {\cdot} A{\cdot}B_m\ ,
\left.\ldots,  \underset{\text{$\ell$-times}}{\underbrace{A{\cdot} \ldots {\cdot}A}}{\cdot}B_1,\ \  \ldots, \ \ \underset{\text{$\ell$-times}}{\underbrace{A{\cdot} \ldots {\cdot}A}}{\cdot}B_m\right \rangle\end{multline}
is equal to $ n =  \dim \cQ$.\par
\end{theo}
\begin{pf} By the discussions in  Examples \ref{example21} and \ref{example44},    the secondary generalised distribution $(V^{II}, \cD^{II})$ of  a system of this kind is (a) regular of rank $n_{\max}$,  (b) involutive and  (c)  $\cM$  admits a decomposition into $\cD^{II}$-strata, consisting of only one stratum, namely  $\cU_0 = \cM$.  Hence, by Theorem \ref{prop93}, any point of $\cM$ is good and the  system is essentially good. \par
Since the vector fields \eqref{voila} are generators of $\cD^{II}$,  any maximal integral leaf of $\cD^{II}$ (which is also  the  $\cD^{II}$-path  component of its points) is an open subset of an affine subspace of $\bR^{1 + n + m}$ of  dimension $m + n_{\max} $ and it is included in a hyperplane $\{ t = \text{const.}\}$. In particular, if $n_{\max} = n$, each such maximal  integral leaf coincides with an  hyperplane 
$\{ t = \text{const.}\}$ and  hence  projects surjectively and homeomorphically onto $\cQ = \bR^n$.  By Corollary \ref{thecor}. the system has the hyper-accessibility property.  \par
Conversely, assume that  $  n_{\text{\rm max}}  <  n$.   By a linear change of coordinates, we may assume that the  subspace  $\bV \subset \bR^n$, which is spanned by the vectors 
\begin{multline*} B_1, \ldots, B_m, \ \ A{\cdot}B_1, \ldots, A{\cdot}B_m, \ \ A {\cdot} A{\cdot}B_1, \ldots, A {\cdot} A{\cdot}B_m, \\
 A {\cdot} A {\cdot} A{\cdot}B_1, \ldots, A {\cdot} A {\cdot} A{\cdot}B_m\ ,
\ldots  \ \end{multline*}
is equal to the subspace $\bV = \{q^{n_{\max}+1} = \ldots, q^n = 0\}$.  Since $\bV$ contains the columns $B_m$  and $  n_{\text{\rm max}} $ is the maximum of the ranks \eqref{seq},   the image of the linear map  $w \mapsto B{\cdot} w$ is  a subspace of $\bV$ and  the linear map $q \mapsto A\cdot q$ leaves   $\bV$ invariant. This implies that  the matrices $A$ and $B$ have the block structures
$$B = \left(\begin{array}{c} \  \ast \  \\
\vdots \\
\  \ast \  \\
 0_{(n - n_{\max}) \times m}
 \end{array} \right)\ ,\  A = \left(\begin{array}{ccc} \ast &  \ast \\
\ast & \ast  \\
0_{(n-n_{\max}) \times n_{\max}} & \wt A
 \end{array} \right)\  \text{for some}\ \wt A \in \bR_{(n-n_{\max}) \times (n-n_{\max})}\ .$$
In particular,  the last $n -n_{\max}$ equations of the system $ \dot q = A q  + B w$ have the form 
\beq \dot q^i = \sum_{j = n_{\max}+1}^n \wt A^i_j q^j\ ,\qquad n_{\max} + 1 \leq i \leq n\ .\eeq
Thus,  for any initial condition $q_o$  with $q_o^{n_{\max}+1} = \ldots = q_o^n =0$ and for any choice of the control curve $w(t)$, the corresponding solution $q(t)$  is a curve which  lies in the subspace  
$\bV = \{q^{n_{\max}+1}, \ldots, q^n = 0\}$.  Hence,  for any such initial conditions, the corresponding reachable set  is included in  $\bV \subsetneq \bR^n$  and   has  empty interior.   This means that the  control system   does not  have  the hyper-accessibility  property. 
\end{pf}
This result  has the  Kalman criterion as a corollary. 
\begin{cor}[Kalman criterion]  \label{corKal} A linear control system  $ \dot q = A q  + B w$, with     $A = (A^i_j) \in \bR_{n \times n}$,  $B = (B^i_a) \in \bR_{n \times m}$ and  $\cK \subset  \bR^m$ is small-time locally controllable at the origin $q_o =  0_{\bR^n}$ if and only if the maximum $n_{\max}$ of the ranks \eqref{seq} is equal to $n = \dim \cQ$. \par
In case $\cK = \bR^m$, the system  is small-time locally controllable  if and only if it has  the global controllability property in the $\cC^\o$ sense (that is, for any $q_o \in \cQ$, the corresponding reachable set  $\bigcup_{T> 0}\Reach^{\cC^\o}_{T}(q_o)$  is the whole state space $\cQ$).
\end{cor}
\begin{pf}  The same  argument of the second part of the proof of Theorem \ref{theo51} shows that if $n_{\max} < n$, there is no solution that joints the  point $q_o = 0_{\bR^n}$  to any point of the complementary set of a proper linear subspace of $\bR^n$.  This implies the necessity of the condition $n_{\max} = \dim \cQ$ in the first  claim.
The sufficiency  comes from Theorem \ref{theo51}  (which guarantees the hyper-accessibility)  and the fact that $q_o = 0_{\bR^n}$ is a state of equilibrium  for the system  $ \dot q = A q  + B w$.\par  For what concerns the second claim, global  controllability trivially implies the small-time local controllability. Conversely,  if the small-time local controllability holds, by the first claim and the proof of Theorem  \ref{theo51}, we have that (a)  $n = n_{\text{max}}$, (b) each maximal integral leaf of the secondary distribution $\cD^{II}$ is a hyperplane $t = \text{const.}$ of $\bR^{1 + n + m}$   and (c) any  maximal surrogate leaflet is an open subset of a hyperplane of this kind.   We now observe that, from   the definition of $\bT^o$-surrogate vector fields and the explicit expression \eqref{voila} of  the generators for the secondary distribution of this control system, any  $\bT^o$-surrogate  field is  a vector field with constant components. In particular its flow is  a one-parameter family of translations of $\bR^{1 + n + m}$, thus  well defined for any  $x \in \cM$ and any value for the parameter $s$ of the flow,  and any pair  of $\bT^o$-surrogate  fields commute. In other words, we have that the   $\bT^o$-surrogate  fields  generate an abelian  group  of translations of dimension $n + m$. Combining this with   the definition of   $\bT^o$-surrogate maps and  of $\bT^o$- surrogate leaflets (Definition \ref{def33}), we immediately get  that each $\bT^o$-surrogate map can be uniquely extended to a map  over the whole $\bR^{n + m}$,  whose image (the corresponding surrogate leaflet)  is the $n+m$-dimensional orbit of the  group of  translations generated by the $\bT^o$-surrogate  fields. This means that any $\bT^o$-surrogate leaflet can be enlarged to become  a whole  hyperplane  $t = \text{const.}$ of $\bR^{1 + n + m}$ and hence to coincide with the $\cD^{II}$-path component of its center.  From the discussion in \S \ref{section32}, we conclude that any $\cM$-attainable set from a point $x_o$ contains all hyperplanes $t = \text{const.}$ passing through its points. Since each such hyperplane  projects surjectively onto $\cQ$,  the global controllability follows. 
\end{pf}
\par
\medskip
\subsection{Three elementary examples}\label{elementary}
Consider the following  very simple non linear control systems, each of them consisting of only one equation: 
\begin{align*}
(A)  \hskip 1 cm  &\dot q = (w^1)^2\\
(B)  \hskip 1 cm  &\dot q = w^1 w^2\\
(C)  \hskip 1 cm &  \dot q = e^{w^1}
\end{align*}
The systems (A) and (C) manifestly do not have the  small-time  local controllability property at any   $q_o$ (by the way, for (A)  any  point $q_o$ is  state of equilibrium). In fact, for $(A)$ and $(C)$, given an initial condition $q_o = q(0)$ and a control curve $w(t)$, the corresponding solution $q(t)$ is  included in the subset  $\{ \ q  \geq q_o\}$  and the reachable set of $q_o$  cannot include any neighbourhood of $q_o$. On the contrary the system (B) does have the small-time local controllability property at any  $q_o$, since   any  reachable set  of $(B)$ is   equal to a   reachable set of the  control system $\dot q = u$ with just one control $u \in \bR$.\par
We want to discuss these three simple systems using Corollary \ref{thecor} and our criterions of goodness. It will offer  a  nice illustration  of  how  our   results work.\par
\smallskip
Let us start  with  $(A)$. In this case, the extended space-time is $\cM = \bR \times \cQ \times \cK$ with $\cQ = \bR$ and $\cK \subset \bR$ and the canonical vector fields are 
\beq \bT^o= \frac{\p}{\p t} + (w^1)^2 \frac{\p}{\p q}\ ,\qquad W^o_1 = \frac{\p}{\p w^1}\ .\eeq
Computing the iterated Lie brackets $[\bT^o, [\bT^o, \ldots,[\bT^o, W^o_1]]$, we find that the set  $V^{II}$  is generated by the   vector fields
\beq W_1^{(0)} = W_1^o = \frac{\p}{\p w^1}\ ,\qquad W_1^{(1)} = [\bT^o, W^o_1] = - 2 w^1 \frac{\p}{\p q}\ .\eeq
The  pair $(W_1^o, W_1^{(1)})$ is  a set of $\bT^o$-adapted generators  with  integers $\nu = 1$, $R_0 =  0$, $R_1 = 1$ and vector fields $W_{\ell(a)j}$, $1 \leq j \leq R_a$, $0 \leq \ell \leq 1$,  equal to  
$$W_{0 (1) 1} \= W^{(0)}_1 = W_1^o\ ,\qquad  W_{1 (1) 1} \= W_1^{(1)} \ .$$
 Looking at these generators, we see that  $\cD^{II}_x$ has dimension $2$ for any $x \in \cM \setminus \{ w^1 = 0\}$ and  has dimension $1$ for  $x \in \{w^1 = 0\}$. Since the  set of vector fields, which is determined by the brackets of vector fields in $V^{II}$, contains not only $W^o_1$, $W_1^{(1)}$,  but also 
$$Y = [W^o_1, W_1^{(1)}] = - 2 \frac{\p}{\p q}\ .$$
 the subspace   $E^{\text{Lie}}_{x} \subset T_x \cM$,  which is spanned by the  values at $x \in \cM$ of the  vector fields in $(V^{II}, \cD^{II})$ and of their iterated Lie brackets,   has dimension $2$ for any $x$.  As it is pointed out in Example \ref{example45}, the whole $\cM$ admits a decomposition into $\cD^{II}$-strata with only one $\cD^{II}$-stratum, i.e.  $\cU_0 = \cM$, and  the corresponding  bracket generated distribution $\cE = \cE^{(\cU_0)}$ is regular and generated by $W^o_1$ and  $Y = [W^o_1, W_1^{(1)}]$.  On the set $\cM \setminus \{w^1 = 0\}$, the generator $Y$ can be replaced by the vector field $W_1^{(1)}$ and this implies that each point of  $\cM \setminus \{w^1 = 0\}$ is a good point of the first kind. But at the points in $\{w^1 = 0\}$ the 
first criterion does not apply and the second criterium neither: Indeed,  there is only one possible sub-distribution to be considered, namely $(V^{II (W^o_1)}, \cD^{II(W^o_1)})$ and this sub-distribution  coincides with $(V^{II}, \cD^{II})$. No point of the set $\{w^1 = 0\}$ is a good point of the first kind and hence  there exists no generator  of the form $Y = [W^{(W^o_\b)} _{B}, W_{B'}]$  satisfying the condition (2) of the second criterion. Therefore,  we may not conclude that the system is essentially good and  Corollary \ref{thecor} cannot be used. \par
As a matter of fact, being the system $(A)$ particularly simple, it is  possible to compute explicitly all possible surrogate leaflets passing through a point of the form $y_o = (t_o, q_o, w_o = 0)$. With a  straightforward (but tedious) computation, one  can directly   check that  no such  point is  good and that our results on surrogate leaflets cannot be used.  This is of course consistent with the above  remarks on the system (A). \par
\medskip
Let us now analyse the  system (B). In this case, the  extended space-time is $\cM = \bR \times \cQ \times \cK$ with $\cQ = \bR$ and $\cK \subset \bR^2$ and the canonical vector fields are 
\beq \bT^o = \frac{\p}{\p t} + w^1 w^2 \frac{\p}{\p q}\ ,\qquad W^o_1 = \frac{\p}{\p w^1}\ , \qquad W^o_2 = \frac{\p}{\p w^2}\ .\eeq
Computing the iterated Lie brackets $[\bT^o, [\bT^o, \ldots,[\bT^o, W_\a^o]]$, $\a = 1,2$,  we see  that   $V^{II}$ is generated by the globally defined vector fields
\beq W^o_1 = \frac{\p}{\p w^1}\ ,\ W_1^{(1)} = [\bT^o, W^o_1] = - 2 w^2 \frac{\p}{\p q}\ ,\ W^o_2 = \frac{\p}{\p w^2}\ ,\ W_2^{(1)} = [\bT^o, W^o_2] = - 2 w^1 \frac{\p}{\p q}\ .\eeq
As for the previous system, also these vector fields form a set of $\bT^o$-adapted generators   for $(V^{II}, \cD^{II})$. The corresponding integers are now $\nu = 1$, $R_0 = 0$,  $R_1 = 2$.\par
\smallskip
Observing that the  generalised distribution,  determined by  the iterated Lie brackets of  the vector fields in $V^{II}$,  contains not only $W^o_\a$, $W_\a^{(1)}$, $\a = 1,2$,  but also 
$$Y = [W^o_1, W_2^{(1)}] = - 2 \frac{\p}{\p q}\ ,$$
the same  arguments for  (A)  lead to the conclusion that also in this case $\cM$ decomposes into $\cD^{II}$-strata with  only one $\cD^{II}$-stratum, i.e. $\cU_0 =  \cM$. We can also directly  see that  all points of $\cM \setminus \{ w^1 = w^2 = 0\}$ are good points of the first kind. On the other hand,  in neighbourhoods of  the points of the set $\{w^1 = w^2 =  0\}$, the bracket generated distribution  contains  the vector fields 
$$W^o_1\ ,\qquad W^o_2\ ,\qquad Y = [W^o_1, W_2^{(1)}]\ .$$
We now remark  that the  vector field $W^o_1$ is in the sub-distribution $(V^{II(W^o_1)}, \cD^{II(W^o_1)})$.  With the usual argument, it is quite simple to check that $(V^{II(W^o_1)}, \cD^{II(W^o_1)})$ is generated by $W^o_1$ and $W_1^{(1)}$,  and that  $\cM$ admits a decomposition into  $\cD^{II(W^o_1)}$-strata with two strata, namely  $\cU_0 = \cM \setminus \{w_2 = 0\}$ and $\cU_1 = \{w_2 = 0\}$. he bracket generated distribution $\cE^{(\cU_0)}$ is regular of rank $2$ and it coincides with $\cD^{II(W^o_1)}|_{\cU_0}$, while  the bracket generated distribution $\cE^{(\cU_1)}$ is regular of rank $1$ and it coincides with $\cD^{II(W^o_1)}|_{\cU_1}$. So, for all points in $\{w^1 = w^2 = 0\}$ (which are all contained in $\cU_1$),  all conditions of the second goodness criterion are satisfied, meaning that  all  these points are good points of the second kind. We therefore conclude that the system (B) is essentially good. By a simple inspection of the maximal  integral leaves of the bracket generated distribution $\cE$ (which are the hyperplanes $\{ t = \text{const.}\}$), it follows that the restriction of $\pi^\cQ$ to any such integral leaf has maximal rank.  Hence  the system has the hyper-accessibility property by Corollary \ref{thecor} and it has   the small-time local controllability property at each   $q_o \in \cQ$,   because any  $q_o \in \cQ$  is a  state of equilibrium. \par
\medskip
Let us now consider the system $(C)$.  It has the same expanded space-time of the system $(A)$, but the canonical vector fields are now
\beq  \bT^o = \frac{\p}{\p t} + e^{w^1} \frac{\p}{\p q}\ ,\qquad W^o_1 = \frac{\p}{\p w^1}\ .\eeq
In this case  the  set $V^{II}$ is  generated by the  vector fields
\beq W^o_1 = \frac{\p}{\p w^1}\ ,\qquad W_1^{(1)} = [\bT^o, W^o_1] = - 2 e^{w^1} \frac{\p}{\p q}\ .\eeq
Moreover the secondary distribution $(V^{II}, \cD^{II})$
is regular of rank $2$ and involutive. For this system, all points are good points of the first kind and the maximal integral leaves of $(V^{II}, \cD^{II})$ are the sets $\{ t = \text{const.} \}$. Thus,  in this case,  Corollary \ref{thecor} applies and the system has the hyper-accessibility property, in contrast with the system $(A)$. However, there is no point $q_o \in \cQ$ that has  the homing property and in particular there is no state of equilibrium. Due to this, we cannot establish the small-time local controllability at any point. This is of course consistent with our initial  observations   on the  system (C). \par
\medskip
\subsection{A locally controllable system  for which the local controllability cannot be established through  the linear test}\label{example3}
 Consider  the     control system with     differential equations 
 \beq \label{marta-1}
\dot q^1=\cos q^3\sin w^1\ ,\qquad 
\dot q^2 =\sin q^3 \sin w^1\ ,\qquad 
\dot q^3=w^2\ , 
\eeq
with controls    $ w(t) =  (w^1(t), w^2(t))$ taking values in an open connected  neighbourhood $\cK$ of $0_{\bR^2}$ in the stripe $\big\{(w^1, w^2)\ :  w^1 \in \big( - \frac{\pi}{2}, \frac{\pi}{2} \big)\big\} \subset \bR^2$. \par
\smallskip
This   system is particularly interesting  for the following reasons.  First,  it  is non-linear in  controls and  thus the  classical Chow controllability criterion (see e.g. \cite[Thm. 3.18]{Co}) is of no use to establish its local controllability. Second, any  $q_o = (q_o^1, q_o^2, q_o^3)$ is a state of equilibrium and, 
for any   corresponding point  $z_o =  (t_o, q^1_o,  q^2_o, q^3_o, 0, 0)$ in the extended space-time $\cM = \bR^4 \times \cK$,  the linearisation at $z_o$ of the  equations \eqref{marta-1}   yields the system of equations on  the translated variables $\wt q^i = q^i - q^i_o$ 
\beq  \frac{d \wt q^1}{dt}  =  \cos q^3_o w^1 \ ,\qquad  \frac{d \wt q^2}{dt} = \sin q^3_o  w^1\ ,\qquad
 \frac{d \wt q^3}{dt} = w^2 \ .
\ \eeq
 According to the  classical Kalman criterion (see e.g. \cite[Thm. 1.16]{Co} or \S \ref{sectKal})  this  linear system (whose solutions approximate  those of the  system \eqref{marta-1} up to  the first order)  has not the controllability property. This means that     the   ``linear test''  (\cite[Thm. 3.8]{Co}) is inconclusive for  the  control  system \eqref{marta-1}.  At the best of our knowledge,  no other   classical  criterium  can be of any help in establishing whether this system  has or does not have  the small-time local controllability property at its points. \par
 \smallskip
We now want to show that our results give an  answer to such a  question. First of all, we observe that for   this example the extended space-time is $\cM = \bR \times \bR^3 \times \cK \subset \bR^6$  and that the canonical vector fields are 
\beq 
\begin{split}
& \bT^o \= \frac{\p}{\p t}  + \cos q^3 \sin w^1 \frac{\p}{\p q^1} +    \sin q^3 \sin w^1\frac{\p}{\p q^2} + w^2 \frac{\p}{\p q^3} \ ,\\
&  W^o_1 \=  \frac{\p}{\p w^1}\ ,\qquad W^o_2 \=  \frac{\p}{\p w^2}\ .
\end{split}
\eeq
Let us now compute   the iterated Lie brackets
$W_\a^{(k)} \=    \underset{k\text{\rm -times}}{\underbrace{[\bT^o, [\bT^o, \ldots, [ \bT^o}}, W^o_\a]\ldots]]$
 in order to determine generators for the secondary distribution $(V^{II}, \cD^{II})$.   In what follows, given a finite set of vector fields $\{X_1, X_2, \ldots, X_r\}$, we denote 
 $$ \big\langle X_1, \ldots , X_r\big\rangle \= \big\{\ X  = f^1 X_1 + \ldots + f^r X_r\  \text{for some real analytic functions} \ f_i\ \big\}\ .$$
 Note that
 $$ \big\{\ X\  \text{or}\  [\bT^o, X]\ \text{for} \ X \in  \langle X_1,  \ldots, X_r \rangle\ \big\} = \big\langle X_1, \ldots, X_r, [\bT^o, X_1],\ldots,  [\bT^o, X_r]\big\rangle\ .$$
 Using this notation, we may write  
 \beq
\begin{split}
&W^{(0)}_1 = W_1^o = \frac{\p}{\p w^1}\ , \\
&W_1^{(1)} =    \cos w^1 \left(\cos q^3 \frac{\p}{\p q^1} - \sin q^3\frac{\p}{\p q^2}\right)\ ,\\
&W_1^{(2)} =  w^2  \cos w^1 \left( \sin q^3 \frac{\p}{\p q^1}  -  \cos q^3 \frac{\p}{\p q^2}\right) \ ,\\
 & W_1^{(3)}  =  0 \hskip -0.2 cm  \mod \langle  W_1^{(1)} \rangle \ ,\\
 \end{split}
 \eeq
  \beq
\begin{split}
&W^{(0)}_2 = W_2^o = \frac{\p}{\p w^2}\ , \\
&W_2^{(1)} =    \frac{\p}{\p q^3} \ ,\\
&W_2^{(2)} =  \sin w^1 \left( \sin q^3 \frac{\p}{\p q^1}  -     \cos q^3 \frac{\p}{\p q^2} \right)\ ,\\
 & W_2^{(3)}  =  - w^2 \sin w^1  \left(\cos q^3 \frac{\p}{\p q^1}  -  \sin q^3 \frac{\p}{\p q^2}\right) \ ,\\
  & W_2^{(4)}  =  0 \hskip -0.2 cm  \mod \langle  W_2^{(2)} \rangle \ .
 \end{split}
 \eeq
 Since $W_1^{(3)}  =  0 \hskip -0.2 cm  \mod \langle  W_1^{(1)} \rangle $ and $W_2^{(4)}  =  0 \hskip -0.2 cm  \mod \langle  W_2^{(2)} \rangle$,  we see that   $(V^{II}, \cD^{II})$ is generated by the globally defined vector fields 
 \beq \label{previous} W^o_1\ ,\ W_1^{(1)}\ ,\ W_1^{(2)}\ ,\ W^o_2\ ,\ W_2^{(1)}\ ,\ W_2^{(2)}\ .\eeq
Since $\cos w^1 \neq 0$ at all points of  $\cM = \bR \times \cQ \times \cK$,  this distribution is regular and of rank $5$  when  it is restricted to $\cM \setminus \text{Sing}$, where
 \beq \text{Sing} = \{\ (t, q^1, q^2, q^2, w^1, w^2) \in \cM\ , \ w^ 1 = w^2 = 0\ \}\ . \eeq
At  any point $x_o \in \text{Sing} $, the space $\cD^{II}_{x_o}$ has the smaller  dimension  $\dim \cD^{II}_{x_o} = 4$.\par
\smallskip 
The bracket generated distribution  $\cE = \cD^{II (\text{Lie})}$ on $\cM$  is regular of rank $5$, being generated by the vector fields in \eqref{previous} together with  the nowhere vanishing vector field
 \beq Y = [ W_1^{(1)}, W_2^{(1)}] =  \cos w^1\left(\sin q^3 \frac{\p}{\p q^1} - \cos q^3 \frac{\p}{\p q^2}\right)\ .\eeq 
 We therefore conclude that  the whole $\cM$  admits a decomposition in $\cD^{II}$-strata, consisting of just  one stratum, namely  $\cU_0 = \cM$. \par
 \smallskip
We now observe that  all points of $ \cM \setminus \text{Sing}$ are good points, because on all sufficiently small neighbourhoods of those points the distribution $\cD^{II}$ is regular. If we can show that also the points of 
 $\text{Sing}$ are good, we may conclude that  the system is essentially good  with  maximal integral leaves of the bracket generated distribution $\cE$ of $(V^{II}, \cD^{II})$, on which the projection $\pi^\cQ$ restricts to a map  of  maximal rank $5 = \dim \cQ$. This would imply that  the system has the hyper-accessibility property by Corollary \ref{thecor} and,   in particular, it has the small time local controllability property at its stable points of equilibrium. \par
 \smallskip
  For this purpose, we observe that the vector field $Y$ is obtained  as a Lie bracket between a vector field in the secondary sub-distribution $(V^{II (W^o_1)}, \cD^{II(W^o_1)})$ and another in the secondary distribution $(V^{II (W^o_2)}, \cD^{II(W^o_2)})$. Moreover, the sub-distribution $(V^{II (W^o_1)}, \cD^{II(W^o_1)})$ admits the  global generators
 \beq W^o_1\ , \ \ W_1^{(1)}\ ,\ \ W_1^{(2)}\ \ .\eeq
 One can directly check that for  this generalised distribution, $\cM$ admits a decomposition into  two $\cD^{II(W^o_1)}$-strata, namely into 
 $$ \cU_0 = \cM \setminus \wt{\text{Sing}} \ ,\qquad \cU_1 = \wt{\text{Sing}} = \{\ (t, q^1, q^2, q^3, w^1, w^2)\in \cM\ \ , \ w^2 = 0\ \}\ .$$
 Indeed, at any point $x_o \in \cU_0$, the space $\cD^{II(W^o_1)}_{x_o}$ has dimension $3$, while at any point  $x_o \in \cU_1$, the space $\cD^{II(W^o_1)}_{x_o}$ has dimension $2$ and it contains only  vectors that are tangent to $\cU_1$.  Moreover, each restriction $\cD^{II(W^o_1)}|_{\cU_i}$, $i = 0, 1$, is not only regular but also involutive.  We thus conclude that   $(V^{II (W^o_1)}, \cD^{II(W^o_1)})$ has the above defined  $\cD^{II(W^o_1)}$-strata and  all conditions of the second criterion are satisfied for any $x_o \in \text{Sing} \subset \wt{\text{Sing}} = \cU_1$. Thus any  such a  point  is good and the system \eqref{marta-1} is essential  good, as we needed  to prove.\par
\medskip
\subsection{Small-time local controllability of a controlled Chaplygin sleigh} \label{example2}
The {\it Chaplygin sleigh} is one of the simplest and interesting examples of nonholonomic dynamical system,  first considered by Chaplygin in \cite{Ch}.  If  external forces and torques are applied, such dynamical system becomes a  particularly interesting nonholonomic  controlled system (see e.g. \cite{BKMM, BK}  and references therein). If the sleigh   is assumed to be   immersed in an ideal  (two-dimensional) potential fluid, the dynamical   system is  called    {\it hydrodynamical Chaplygin sleigh} and, in first approximation, it  is  controlled  according to a    system of   equations of the form
 \beq \label{sleigh-1}
\begin{split}
&\dot q^1=q^4\cos q^3\ ,\qquad 
\dot q^2 =q^4\sin q^3\ ,\qquad 
\dot q^3=q^5\ ,\\
&\dot q^4=w^1+ {\mathbb A}\, q^4 q^5\ ,\qquad 
\dot q^5=w^2\ ,
\end{split}
\eeq
where    $\mathbb A$ is a constant, $(q^1, q^2, q^3, q^4, q^5)$ are state coordinates for the system (corresponding to the position, the center of mass speed and the angular velocity of the sleigh) and   $ w =  (w^1, w^2)$ are control parameters, corresponding to external forces  and torques acting on the dynamical system.  The controls $w$  are assumed to  take values in an open connected subset $\cK$ of  $\bR^2$. The constant $\mathbb A$ is what distinguishes a hydrodynamical sleigh from the classical  sleigh: When $\bA = 0$, the equations reduce to those of the classical controlled Chaplygin sleigh.   For details on the physical meanings of the state and  coordinates   of this dynamical system and for   investigations on various aspects, we refer to  \cite{Ba, FG,  FGV, SZ} and references therein.\par
\smallskip
To the best of our knowledge, besides research efforts of various type, the
(global) controllability or  the small-time local controllability of a dynamical system,   guided by the system of equations \eqref{sleigh-1},   has not been proved so far,  not even
for the case in which $\bA = 0$, i.e. for  the  classical controlled Chaplygin sleigh. The goal of this section is to show that such dynamical system is indeed small-time locally controllable  at its states of equilibrium by virtue of our  criterions for goodness. Studies on the global controllability are left to future investigations.\par
\smallskip 
First of all, we observe that for a dynamical system controlled by \eqref{sleigh-1},  the extended space-time is $\cM = \bR \times \bR^5 \times \cK \subset \bR^8$  and the associated canonical vector fields are  
\beq \label{bt-1}
\begin{split}
& \bT^o \= \frac{\p}{\p t}  + q^4 \cos q^3 \frac{\p}{\p q^1} +   q^4 \sin q^3 \frac{\p}{\p q^2} + q^5 \frac{\p}{\p q^3} 
+\left(w^1 +  {\mathbb A}\, q^4 q^5\right) \frac{\p}{\p q^4} + w^2 \frac{\p}{\p q^5}\ ,\\
&W^o_1 \= \frac{\p}{\p w^1} \ ,\qquad  W^o_2 \=  \frac{\p}{\p w^2}\ \ .\end{split}
\eeq
Let us compute   the  Lie brackets
$W_\a^{(k)} \=\underset{k-\text{times}}{\underbrace{ [\bT^o, [\bT^o, \ldots, [ \bT^o, }}W^o_\a]\ldots]]$
 and  determine generators for the secondary distribution $(V^{II}, \cD^{II})$.  As in \S \ref{example3}, given a  set of vector fields $\{X_i\}_{1 \leq i \leq r}$, we  set 
 $ \big\langle X_1, \ldots , X_r\big\rangle \= \big\{\ X  = \sum_{i = 1}^r f^i X_i\ \text{for some $\cC^\o$ functions} \ f_i \ \big\}$.   We get
\beq
\nonumber
\begin{split}
&W^o_1 = \frac{\p}{\p w^1}\ , \\
&W_1^{(1)} =  - \frac{\p}{\p q^4}\ ,\\
&W_1^{(2)} =   \cos q^3 \frac{\p}{\p q^1} +    \sin q^3 \frac{\p}{\p q^2} +\bA q^5 \frac{\p}{\p q^4}  =   \\
& \hskip 5 cm  = \cos q^3 \frac{\p}{\p q^1} +    \sin q^3 \frac{\p}{\p q^2} \hskip -0.2 cm  \mod \langle W_1^{(1)} \rangle\\
 & W_1^{(3)}  =  q^5\left(- \sin q^3 \frac{\p}{\p q^1} + \cos q^3 \frac{\p}{\p q^2} \right) \hskip -0.2 cm  \mod \langle 
  W_1^{(1)}, W_1^{(2)} \rangle\\
 &  W_1^{(4)} =  w^2\left(- \sin q^3 \frac{\p}{\p q^1} + \cos q^3 \frac{\p}{\p q^2} \right) \hskip -0.2 cm  \mod \langle  W_1^{(1)}, W_1^{(2)}, W^{(3)}_1\rangle\ ,\\
  &  W_1^{(5)} =  0 \hskip -0.2 cm  \mod \langle  W_1^{(1)}, W_1^{(2)}, W^{(3)}_1, W^{(4)}_1 \rangle\ ,\\
&W^o_2 = \frac{\p}{\p w^2}\ , \\
&W_2^{(1)} =  - \frac{\p}{\p q^5}\ ,\\
&W_2^{(2)} =    \frac{\p}{\p q^3} + \bA q^4 \frac{\p}{\p q^4} = \frac{\p}{\p q^3} \hskip -0.2 cm \mod \langle  W_1^{(1)}\rangle\ ,\\
 &W_2^{(3)}   = q^4\left(\sin q^3 \frac{\p}{\p q^1} - \cos q^3 \frac{\p}{\p q^2}\right) \hskip -0.2 cm  \mod \langle  W_1^{(1)}, W_1^{(2)}\rangle\ ,\\
 & W^{(4)}_2  = q^4 q^5 \left( \cos q^3 \frac{\p}{\p q^1} + \sin q^3 \frac{\p}{\p q^2}\right) +\\
 & \hskip 1 cm  + \left( w^1 + \bA q^4 q^5\right) \left( \sin q^3 \frac{\p}{\p q^1} - \cos q^3 \frac{\p}{\p q^2}\right)  \hskip -0.2 cm  \mod \langle  W_1^{(1)}, W_1^{(2)}, W_1^{(3)}\rangle = 
 \\
 &  \hskip 2 cm =  w^1  \left( \sin q^3 \frac{\p}{\p q^1} - \cos q^3 \frac{\p}{\p q^2}\right)  \hskip -0.2 cm  \mod \langle  W_1^{(1)}, W_1^{(2)}, W_1^{(3)}\rangle\ , \\
 & W^{(5)}_2 =  w^1 q^5 \left(\cos q^3 \frac{\p}{\p q^1} +  \sin q^3 \frac{\p}{\p q^2} \right)  \hskip -0.2 cm  \mod \langle  W_1^{(1)}, W_1^{(2)}, W_1^{(3)}, W_1^{(4)}\rangle =  \\
  & \hskip 3 cm = 0\hskip -0.2 cm  \mod \langle  W_1^{(1)}, W_1^{(2)}, W_1^{(3)}, W_1^{(4)}\rangle \ .
\end{split}
\eeq
 We  see that  that $(V^{II}, \cD^{II})$ is globally generated by the vector fields 
\beq \label{genchap} W^o_1\ ,\ \  W_1^{(1)}\ ,\ \  W_1^{(2)}\ ,\ \  W_1^{(3)}\ ,\ \  W_1^{(4)}\ ,\ \ W^o_2\ ,\ \ W_2^{(1)}\ ,\ \  W_2^{(2)}\ ,\ \ W_2^{(3)}\ ,\ \  W_2^{(4)}\ .\eeq
Note that, for any  $y_o \in \cM \setminus \{q_5 \neq 0\}$, the  set  $\{W^o_1\ ,\  W_1^{(1)}\ ,\  W_2^{(2)}\ ,\ W_1^{(3)}\ ,\ W^o_2\ ,\ W_2^{(1)}\ ,\  W_2^{(2)}\}$ is  a set of  $\bT^o$-adapted generators on a neighbourhood $\cU$ of $y_o$ 
with corresponding integers $\nu = 3$,     $R_0 = R_1 = 0$ and $R_2 = R_3 = 1$.  Other choices of generators among the  \eqref{genchap}   provide $\bT^o$-adapted generators on neighbourhoods of    points  of  $\{q^5 = 0\}$.  \par
Denoting by $\text{Sing} = \{\ x \in \cM\ : \ q^4 = q^5 = w^1 = w^2 = 0\}$, we have that 
\begin{itemize}
\item[--] if $x \in \text{Sing}$, then $\dim \cD^{II}_x = 6$; 
\item[--] if $x \in \cM \setminus \text{Sing}$, then $\dim \cD^{II}_x = 7 = \dim \cQ \times \cK$;
\end{itemize}
On the other hand, in any set  of vector fields   containing all  Lie brackets between  vector fields in the set  \eqref{genchap} there is the vector field 
\begin{multline}  Y = [W_1^{(3)}, W_2^{(1)}] = \sin q^3 \frac{\p}{\p q^1} - \cos q^3 \frac{\p}{\p q^2}\\
 \hskip - 0.2 cm \mod\big \langle W_1^{(1)}, W_1^{(2)}, [W_1^{(1)}, W_2^{(1)}] = 0, [W_1^{(2)}, W_2^{(1)}] = 0\big\rangle\ .\end{multline}
 Since  $Y$ and  the vector fields in \eqref{genchap} generate a regular distribution of rank $7$,   given $x_o \in \cM$ (regardless on whether it is in \text{Sing} or not), the  corresponding space $\cE^{\text{Lie}}_{x_o} \subset T_{x_o} \cM$, generated by  all Lie brackets of vector fields in $V^{II}$,   has constantly  dimension $7$.  Hence, the whole $\cM$ decomposes into a single $\cD^{II}$-stratum,  namely into $\cU_0 = \cM$,  and the corresponding bracket generated distribution $\cE$ has  constant rank $7$. The maximal integral leaves of $\cE$ are the intersections of $\cM = \bR^6 \times \cK \subset \bR^8 $ with the  hyperplanes $\{\ t = \text{const.}\ \}$.\par
 All points of $ \cM \setminus \text{Sing}$ are good points of the first kind, because on  sufficiently small neighbourhoods of those points the distribution $\cD^{II}$ is regular. We claim that  the points of $\text{Sing}$ are good of the second kind and hence that the dynamical system is essentially good. To check this, we first observe that the vector field $Y$ is  a Lie bracket between a vector field in the secondary sub-distribution $(V^{II (W^o_1)}, \cD^{II(W^o_1)})$ and one in the secondary distribution $(V^{II (W^o_2)}, \cD^{II(W^o_2)})$.  The sub-distribution $(V^{II (W^o_1)}, \cD^{II(W^o_1)})$ admits the  set of (global defined) generators
 \beq W^o_1\ , \ \ W_1^{(1)}\ ,\ \ W_1^{(2)}\ ,\ \ W_1^{(3)}\ ,\ \ W_1^{(4)}\ .\eeq
 At any point   $x \in \cM \setminus \{ q^5 = w^2 = 0\}$, the space $\cD^{II(W^o_1)}_x$ has dimension $5$, while at the points of $\{q^5 = w^2 = 0\}$, the space $\cD^{II(W^o_1)}_x$ has dimension $4$.  Note also that, for any  $x_o \in \{q^5 = w^2 = 0\}$, the corresponding  space $\cE_{x_o}^{\text{Lie}} \subset T_{x_o} \cM$, generated by the iterated Lie brackets of vector fields in $\cD^{II(W^o_1)}_x$, is  tangent to the hyperplane $\{q^5 = w^2 = 0\}$.  Hence there are two good candidates for the $\cD^{II(W^o_1)}$-strata, namely $\cU_0 = \cM \setminus \{q^5 = w^2 = 0\}$ and $\cU_1 =\{q^5 = w^2 = 0\}$. Computing the    corresponding bracket generated distributions $\cE^{(\cU_i)}$, $ i = 0, 1$, one can easily see that they are both regular  on the corresponding strata. We therefore conclude that   $(V^{II (W^o_1)}, \cD^{II(W^o_1)})$ is of stratified uniform type with the above defined  $\cD^{II(W^o_1)}$-strata. One can also check that for any  $x_o\in \text{Sing} \subset \{q^5 = w^2 = 0\}$ all other conditions of the second criterion are satisfied, meaning that   $x_o$ is a good point of the second kind, as claimed. 
It remains now to check that  the rank of  the restriction of the projection $\pi^\cQ$ on any maximal integral leaves of the bracket generated distribution $\cE$ of $(V^{II}, \cD^{II})$ is maximal.  A direct inspection shows that it is  equal to $5 = \dim \cQ$ at any point. The system has therefore the hyper-accessibility property by Corollary \ref{thecor} and, consequently,  it is  small-time locally controllable  at its states of equilibrium.\par
\ \\[-37pt]

\ \\[-9pt]
\font\smallsmc = cmcsc8
\font\smalltt = cmtt8
\font\smallit = cmti8
\!\!\hbox{\parindent=0pt\parskip=0pt
\vbox{\baselineskip 9.5 pt \hsize=3.7truein
\obeylines
{\smallsmc 
Marta Zoppello
Dipartimento di Scienze Matematiche 
``G. L. Lagrange'' (DISMA)
Politecnico di Torino
Corso Duca degli Abruzzi, 24, 
10129 Torino 
ITALY}\medskip
{\smallit E-mail}\/: {\smalltt  marta.zoppello@polito.it 
\ 
}
}
\hskip -0.8truecm
\vbox{\baselineskip 9.5 pt \hsize=3.7truein
\obeylines
{\smallsmc
Cristina Giannotti \& Andrea Spiro
Scuola di Scienze e Tecnologie
Universit\`a di Camerino
Via Madonna delle Carceri
I-62032 Camerino (Macerata)
ITALY
\ 
}\medskip
{\smallit E-mail}\/: {\smalltt cristina.giannotti@unicam.it
\smallit E-mail}\/: {\smalltt andrea.spiro@unicam.it}
}
}

\begin{thebibliography}{25}
%
%
%
%
%
 \bibitem{AS}
A.\ A. Agrachev and  Y.\ L. Sachkov,  {Control theory from the geometric viewpoint} in ``Control Theory and Optimization  II'', {\it Springer-Verlag, Berlin}, 2004.
%
\bibitem{Ba} C.\ Barot, The motion and control of a Chaplygin Sleigh with internal
shape in an ideal fluid, {\it Ph.D. dissertation, Elect. Eng., Univ. Michigan,
Ann Arbor, MI, USA}, 2017.
%
              
%
  
\bibitem{Bl}
A. Bloch,  Nonholonomic mechanics and control, 
 {\it Springer, New York},  2015.
 
\bibitem{BKMM}  A.\ M.\ Bloch, P.\ S.\ Krishnaprasad, J.\ E.\ Marsden, and R. Murray,
{\it Nonholonomic mechanical systems with symmetry}, Arch. Ration. Mech. Anal.  {\bf 136} (1996), 21--99.
 
 
 \bibitem{BK} 
 A.\ V. \ Borisov and S.\ P. \ Kuznetsov, {\it Regular and chaotic motions of
a Chaplygin Sleigh under periodic pulsed torque impacts},  Regul.
Chaotic Dyn.  {\bf 21} (2017), 792--803.

   
\bibitem{BR} A. Bressan and F. Rampazzo, {\it On differential systems with vector-valued impulsive controls},
Boll. Un. Mat. Ital. B (7) {\bf 2}, (1988),  641--656. 
            
 
            
              
 

 
\bibitem{Ch} S.\ A.\ Chaplygin, {\it On the theory of motion of nonholonomic systems.
The reducing-multiplier theorem} (Russian),   Mat. Sb. {\bf 28} (1911) (English translation  in Regul. Chaotic Dyn. {\bf 13} (2008)
369--376).


\bibitem{Co}
J.\ M.\ Coron,  Control and nonlinearity, {\it  American Mathematical Society, Providence, RI}, 2007.
 
 \bibitem{FG} Y.\ N.\ Fedorov and L.\ CX.\  García-Naranjo, {\it The hydrodynamic Chaplygin
Sleigh}, J. Phys. A Math. Theor. {\bf 43} (2010) Art. no. 434013.
 
\bibitem{FGV}   Y.\ N.\ Fedorov,  L.\ C.\ Garc\'{\i}a-Naranjo and J. 
              Vankerschaver,
 {\it The motion of the 2{D} hydrodynamic {C}haplygin sleigh in the
              presence of circulation},
 Discrete Contin. Dyn. Syst.  {\bf 33} (2013),
4017--4040.
      
      
 

  
  \bibitem{GSZ1} C. Giannotti, A. Spiro and M. Zoppello, {\it  Proving the  Chow-Rashevski\u\i\ Theorem \`a la Rashevski\u\i}, in Arxiv 2401.07546 (2024).
     \bibitem{GSZ2} C. Giannotti, A. Spiro and M. Zoppello, {\it Distributions and controllability problems (I)}, in Arxiv 2401.07555  (2024). 


     
%
%
%
 
\bibitem{Ma} B. Malgrange, 
 Ideals of differentiable functions,
   {\it Tata Institute of Fundamental Research, Bombay; Oxford
              University Press, London}, 1967.
     
 
 
\bibitem{Na} T. Nagano,
{\it Linear differential systems with singularities and an
              application to transitive {L}ie algebras},
J. Math. Soc. Japan  {\bf 18},
(1966),
  398--404.
 
 
 
\bibitem{OZ}  J.\ M.\ Osborne and D.\ V.\  Zenkov, {\it Steering the Chaplygin Sleigh by a
moving mass} in ``Proceedings of the 44th IEEE Conference on Decision and Control'', Seville, Spain, 2005, pp. 1114--1118,  {\it  Proc. CDC}, 2005.



 \bibitem{Ra}
P. K. Rashevski\u\i, {\it About connecting two points of complete non-holonomic space by admissible curve} (in Russian), Uch. Zapiski Ped. Inst. K. {\bf 
2} 
(1938), 83 -- 94.
 
 
 
 \bibitem{SZ}
N. Sansonetto and  M.  Zoppello,   {\it On the trajectory generation of the hydrodynamic {C}haplygin sleigh}, IEEE Control Syst. Lett.  {\bf 4} (2020), 922--927.
  
  \bibitem{Su} H.\ J.\ Sussmann, {\it Orbits of families of vector fields and integrability of distributions}, Trans. Amer. Math. Soc. {\bf 180}
(1973), 171--188.


\bibitem{SJ} H.\ J.\ Sussmann and V.\ Jurdjevic,
    {\it Controllability of nonlinear systems},  J. Differential Equations  {\bf 12},
(1972), 95--116.


\bibitem{Wa} F.\ W.\ Warner, Foundations of Differentiable Manifolds and Lie groups, {\it Springer-Verlag, New York}, 1983.

 \bibitem{Ju} V.\ Jurdjevic, Geometric control theory,
 {\it Cambridge University Press, Cambridge}, 1997.
%
%
%
%
\end{thebibliography}
\end{document}